\documentclass[12pt]{amsart}

\usepackage[margin=1.0in]{geometry}
\usepackage[english]{babel}
\usepackage[utf8]{inputenc}
\usepackage{fancyhdr}
\usepackage{natbib}
\usepackage{amsmath}

\RequirePackage{amsthm,amsmath,amsfonts,amssymb}
\usepackage{graphicx}
\usepackage{hyperref}
\usepackage{graphicx}
\usepackage{placeins}
\usepackage{caption}
\usepackage{subcaption}
\numberwithin{equation}{section}

\theoremstyle{plain}
\theoremstyle{remark} 
\numberwithin{equation}{section} 

\title{Some Unified Results on Isotonic Regression Estimators of Order Restricted Parameters of a General Bivariate Location/Scale Model}

\author{Naresh Garg  and Neeraj Misra \\ {\footnotesize Department of Mathematics and Statistics\\Indian Institute of Technology Kanpur \\Kanpur-208016, Uttar Pradesh, India}}

\makeatletter
\def\@seccntformat#1{%
  \protect\textup{\protect\@secnumfont
    \ifnum\pdfstrcmp{subsection}{#1}=0 \bfseries\fi
    \csname the#1\endcsname
    \protect\@secnumpunct
  }%
}  
\makeatother

\begin{document}
\maketitle
\section*{\textbf{Abstract}}

	We consider component-wise estimation of order restricted location/scale parameters $\theta_1$ and $\theta_2$ ($\theta_1\leq \theta_2$) of a general bivariate distribution under the squared error loss function. To find improvements over the best location/scale equivariant estimators (BLEE/BSEE) of $\theta_1$ and $\theta_2$, we study isotonic regression of suitable location/scale equivariant estimators (LEE/SEE) of $\theta_1$ and $\theta_2$ with general weights. Let $\mathcal{D}_{1,\nu}$ and $\mathcal{D}_{2,\beta}$ denote suitable classes of isotonic regression estimators of $\theta_1$ and $\theta_2$, respectively. Under the squared error loss function, we characterize admissible estimators within classes $\mathcal{D}_{1,\nu}$ and $\mathcal{D}_{2,\beta}$, and identify estimators that dominate the BLEE/BSEE of $\theta_1$ and $\theta_2$. Our study unifies and extends several studies reported in the literature for specific probability distributions having independent marginals. Additionally, some new and interesting results are obtained. A simulation study is also considered to compare the risk performances of various estimators.
\\~\\ \textbf{Keywords:} Admissible estimators; BLEE; BSEE; Inadmissible estimators; Isotonic regression; Location/Scale model; Mixed estimators \and Scaled squared error loss (SSEL); Squared error loss (SEL).



\section{\textbf{Introduction}}\label{intro}
Let $\underline{X}=(X_1,X_2)$ be a random vector having the Lebesgue probability density function (pdf) $f_{\underline{\theta}}(x_1,x_2)$, $(x_1,x_2)\in \Re^2$, where $\underline{\theta}=(\theta_1,\theta_2)\in \Theta\subseteq \Re^2$ is a vector of unknown parameters and $\Theta$ is the parameter space (to be called unrestricted parameter space); here $\Re$ denotes the real line and $\Re^2=\Re\times \Re$. For a completely specified pdf $f(\cdot,\cdot)$, we assume one of the following models for $f_{\underline{\theta}}(x_1,x_2)$:
\begin{align}
	f_{\underline{\theta}}(x_1,x_2)&= 	f(x_1-\theta_1,x_2-\theta_2),\;(x_1,x_2)\in \Re^2, \;\underline{\theta}=(\theta_1,\theta_2)\in \Theta=\Re^2,\label{eq:1.1}
	\\
	\text{or}\quad	f_{\underline{\theta}}(x_1,x_2)&=	\frac{1}{\theta_1 \theta_2}f\!\left(\frac{x_1}{\theta_1},\frac{x_2}{\theta_2}\right),\; (x_1,x_2)\in \Re^2, \;\underline{\theta}=(\theta_1,\theta_2)\in \Theta=\Re_{++}^2,\;\;\,\label{eq:1.2}
\end{align} 

\noindent 
where $\Re_{++}=(0,\infty)$, $\Re_{++}^2=\Re_{++}\times \Re_{++}$, and $\Theta$ denotes the parameter space. The models \eqref{eq:1.1} and \eqref{eq:1.2} will be referred to as the location and scale probability models, respectively, and the corresponding parameters $\theta_1$ and $\theta_2$ will be called the location and the scale parameters, receptively. Generally, $\underline{X}=(X_1,X_2)$ would be a minimal sufficient statistic based on a random sample from a bivariate distribution or a minimal sufficient statistic based  on independent random samples. Under the scale model \eqref{eq:1.2}, we will assume, throughout, that the distributional support of $\underline{X}=(X_1,X_2)$ is a subset of $\Re_{++}^2$.\vspace*{2mm}

Consider estimation of $\theta_i$ under the loss function $L_i(\underline{\theta},a),\;\underline{\theta}\in\Theta,\;a\in \mathcal{A},\;i=1,2,$ where $\mathcal{A}$ denotes the action space. While working under the location model \eqref{eq:1.1}, we assume that
\begin{equation}\label{eq:1.3}
	L_i(\underline{\theta},a)=(a-\theta_i)^2,\;\; \;\underline{\theta}\in\Theta,\;a\in\mathcal{A}=\Re,\;i=1,2,
\end{equation}
and, under the scale model \eqref{eq:1.2}, we take
\begin{equation}\label{eq:1.4}
	L_i(\underline{\theta},a)=\left(\frac{a}{\theta_i}-1\right)^2,\;\; \;\underline{\theta}\in\Theta,\;a\in\mathcal{A}=\Re_{++},\;i=1,2.
\end{equation}
The loss functions \eqref{eq:1.3} and \eqref{eq:1.4} will be referred to as the squared error loss (SEL) function and the scaled squared error loss (SSEL) function, respectively.
\vspace*{2mm}

The problem of estimating location parameter $\theta_i$, under the location model \eqref{eq:1.1} and loss function \eqref{eq:1.3}, is invariant under the additive group of transformations $\mathcal{G}_L=\{g_{c_1,c_2}:(c_1,c_2)\in\Re^2 \}$, where $g_{c_1,c_2}(x_1,x_2)=(x_1+c_1,x_2+c_2)$, $(x_1,x_2)\in\Re^2,\;  (c_1,c_2)\in\Re^2$. Any location equivariant estimator (LEE) of $\theta_i$ has the form $\delta_{i,c}(\underline{X})=X_i-c,\;c\in\Re,\;i=1,2$. The best location equivariant estimator (BLEE) of $\theta_i$ is $\delta_{i,c_{0,i}}(\underline{X})=X_i-c_{0,i}$, where 
\begin{equation}\label{eq:1.5}
	c_{0,i}\!=E[Z_i]=\!\int_{-\infty}^{\infty}\int_{-\infty}^{\infty}\!\!z_if(z_1,z_2)dz_1 dz_2,
\end{equation}
for $Z_i=X_i-\theta_i,\;i=1,2$. Similarly, estimation of $\theta_i$, under the scale model \eqref{eq:1.2} and the SSEL function \eqref{eq:1.4}, is invariant under the multiplicative group of transformations $\mathcal{G}_S=\{g_{c_1,c_2}:(c_1,c_2)\in\Re_{++}^2\}$, where $g_{c_1,c_2}(x_1,x_2)=(c_1x_1,c_2x_2)$, $(x_1,x_2)\in\Re^2,\; (c_1,c_2)\in\Re_{++}^2$. Any scale equivariant estimator (SEE) of $\theta_i$ is of the form $\delta_{i,c}(\underline{X})=c X_i,\;c>0,\;i=1,2$. The best scale equivariant estimator (BSEE) of $\theta_i$ is $\delta_{i,c_{0,i}}(\underline{X})=c_{0,i} X_i$, where, for $Z_i=\frac{X_i}{\theta_i}$,
\begin{equation}\label{eq:1.6}
	c_{0,i}\!=\frac{E[Z_i]}{E[Z_i^2]}=\!\frac{\int_{-\infty}^{\infty}\int_{-\infty}^{\infty}z_if(z_1,z_2)dz_1 dz_2}{\int_{-\infty}^{\infty}\int_{-\infty}^{\infty}z_i^2f(z_1,z_2)dz_1 dz_2},\;i=1,2.
\end{equation}
Under the unrestricted parameter space $\Theta$, the BLEE and BSEE of $\theta_i$ are known to possess various optimality properties.
\vspace*{1.5mm}

In many situations, it may be known apriori that  $\theta_1\leq \theta_2$ so that the parameter space of interest is $\Theta_0=\{(x_1,x_2)\in\Theta:x_1\leq x_2\}$, called the restricted parameter space. In such situations, it may be of interest to estimate $\theta_1$ and $\theta_2$. For descriptions of various real-life situations where such an estimation problem may be of interest, one may refer to Barlow et. al. (\citeyear{MR0326887}), Robertson et. al. (\citeyear{MR961262}) and Kumar and Sharma (\citeyear{MR981031}). In such problems, a natural question that arises is whether the prior information $\theta_1\leq \theta_2$ can be exploited to find better alternatives to the BLEE/BSEE. This issue has been extensively investigated in the literature for specific probability distributions, having independent marginals (see, for example,  Kumar and Sharma (\citeyear{MR981031}), Kushary and Cohen (\citeyear{MR1029476}), Pal and Kushary (\citeyear{MR1165709}), Kubokawa and Saleh (\citeyear{MR1370413}), Misra and Singh (\citeyear{MR1366828}) and Chang et al. (\citeyear{MR3698503})). An account of these studies can be found in van Eeden (\citeyear{MR2265239}). \vspace*{2mm}

To find improvements over the BLEE/BSEE, we will consider isotonic regression of suitable pair of LEE/SEE of $\underline{\theta}=(\theta_1,\theta_2)$. Let $\omega_1$ and $\omega_2$ be given non-negative weights, such that $\omega_1+\omega_2\neq 0$, and let $(\delta_{1,\beta}(\underline{X}),\delta_{2,\nu}(\underline{X}))$ be a given pair of unrestricted LEE/SEE of $\underline{\theta}=(\theta_1,\theta_2)$ ($\delta_{1,\beta}(\underline{X})\leq\delta_{2,\nu}(\underline{X})$, almost surely, does not hold). The isotonic regression of $(\delta_{1,\beta}(\underline{X}),\delta_{2,\nu}(\underline{X}))$, with weights $(\omega_1,\omega_2)$, is defined as the pair $(\delta^{*}_{1,\beta}(\underline{X}),\delta^{*}_{2,\nu}(\underline{X}))$ that minimizes, for each sample point $\underline{x},$
$$Q=\omega_1 (\delta_{1,\beta}(\underline{x})-y_1(\underline{x}))^2+\omega_2 (\delta_{2,\nu}(\underline{x})-y_2(\underline{x}))^2,$$
among all $y_1(\underline{x})$ and $y_2(\underline{x})$ such that $y_1(\underline{x})\leq y_1(\underline{x})$. For $\alpha=\frac{\omega_1}{\omega_1+\omega_2}$, it turns out that (see Barlow et al. (1972))
\begin{align*}	\delta^{*}_{1,\beta}(\underline{X})&=\min\{\delta_{1,\beta}(\underline{X}),\alpha \delta_{1,\beta}(\underline{X})+(1-\alpha)\delta_{2,\nu}(\underline{X})\}\\
	&=\begin{cases}
		\delta_{1,\beta}(\underline{X}), &\text{ if}\quad \delta_{1,\beta}(\underline{X})\leq \delta_{2,\nu}(\underline{X}) \\
		\alpha \delta_{1,\beta}(\underline{X})+(1-\alpha)\delta_{2,\nu}(\underline{X}),&\text{ if}\quad \delta_{1,\beta}(\underline{X})> \delta_{2,\nu}(\underline{X})
	\end{cases} 
\end{align*}
and                                                                                    \begin{align*}	\delta^{*}_{2,\nu}(\underline{X})&=\max\{\delta_{2,\nu}(\underline{X}),\alpha \delta_{1,\beta}(\underline{X})+(1-\alpha)\delta_{2,\nu}(\underline{X})\}\\
	&=\begin{cases}
		\delta_{2,\nu}(\underline{X}), &\text{ if}\quad \delta_{1,\beta}(\underline{X})\leq \delta_{2,\nu}(\underline{X}) \\
		\alpha \delta_{1,\beta}(\underline{X})+(1-\alpha)\delta_{2,\nu}(\underline{X}),&\text{ if}\quad \delta_{1,\beta}(\underline{X})> \delta_{2,\nu}(\underline{X})
	\end{cases} .
\end{align*}                                                                             We call the above estimators as isotonic regression estimators based on $(\delta_{1,\beta},\delta_{2,\nu})$. For a detailed treatment of the notion of isotonic regression, readers are referred to Barlow et. al. (\citeyear{MR0326887}) and Robertson et. al. (\citeyear{MR961262}).\vspace*{2mm}

Let $c_{0,1}$ and $c_{0,2}$ be as defined by \eqref{eq:1.5}/\eqref{eq:1.6}, so that $\delta_{i,c_{0,i}}$ is the BLEE/BSEE of $\theta_i,\;i=1,2$. To find improvements over the BLEE/BSEE $\delta_{1,c_{0,1}}$ of $\theta_1$, under the prior information $\underline{\theta}\in \Theta_0$, we propose to consider isotonic regression $(\delta_{1,c_{0,1},\alpha}^{(\nu)},\delta_{2,c_{0,1},\alpha}^{(\nu)})$ of $(\delta_{1,c_{0,1}},\delta_{2,\nu})$, for appropriately chosen $\nu$ and weight $\alpha$, and take $\delta_{1,c_{0,1},\alpha}^{(\nu)}$ as a competitor of the BLEE/BSEE $\delta_{1,c_{0,1}}$. Similarly, for estimation of $\theta_2$, under the prior information $\underline{\theta}\in \Theta_0$, we propose to consider isotonic regression $(\delta_{1,c_{0,2},\alpha}^{(\beta)},\delta_{2,c_{0,2},\alpha}^{(\beta)})$ of $(\delta_{1,\beta},\delta_{2,c_{0,2}})$, for appropriately chosen $\beta$ and weight $\alpha$, and take $\delta_{2,c_{0,2},\alpha}^{(\beta)}$ as a competitor of the BLEE/BSEE $\delta_{2,c_{0,2}}$. Here, one may argue that natural choices for $\nu$ and $\beta$ are $c_{0,2}$ and $c_{0,1}$ (the one corresponding to the BLEE/BSEE of $\theta_2$ and $\theta_1$, respectively). As we will observe later that such a choice of $\nu$ or $\beta$ may not always be successful in providing improvements over the BLEE/BSEE $\delta_{1,c_{0,1}}$ or $\delta_{2,c_{0,2}}$, as desired (see examples provided in Sections 3.3 and 4.3). For this reason, to begin with, we take $\nu$ or $\beta$ to be a fixed constant and, depending on the application at hand, we choose them appropriately to achieve our goal of improving the BLEE/BSEE. The guiding principle for choosing $\nu$ or $\beta$ in any application is that $\delta_{2,\nu}$ or $\delta_{1,\beta}$ are close to $\delta_{2,c_{0,2}}$ or $\delta_{1,c_{0,1}}$ and the required conditions (as derived in various results proved in the paper) for improving on the BLEE/BSEE are satisfied. This leads to studying isotonic regression estimators
\begin{align}\label{eq:1.7}
	\delta_{1,c_{0,1},\alpha}^{(\nu)}(\underline{X})=\begin{cases}
		\delta_{1,c_{0,1}}(\underline{X}), &\text{if}\; \delta_{1,c_{0,1}}(\underline{X})\leq \delta_{2,\nu}(\underline{X}) \\
		\alpha \delta_{1,c_{0,1}}(\underline{X})+(1-\alpha)\delta_{2,\nu}(\underline{X}),&\text{if}\; \delta_{1,c_{0,1}}(\underline{X})> \delta_{2,\nu}(\underline{X})
	\end{cases} , \alpha\in[0,1],
\end{align}
for estimating $\theta_1$, and isotonic regression estimators 
\begin{align}\label{eq:1.8}
	\delta_{2,c_{0,2},\alpha}^{(\beta)}(\underline{X})=\begin{cases}
		\delta_{2,c_{0,2}}(\underline{X}), &\text{if}\; \delta_{1,\beta}(\underline{X})\leq \delta_{2,c_{0,2}}(\underline{X}) \\
		\alpha \delta_{1,\beta}(\underline{X})+(1-\alpha)\delta_{2,c_{0,2}}(\underline{X}),&\text{if}\; \delta_{1,\beta}(\underline{X})> \delta_{2,c_{0,2}}(\underline{X})
	\end{cases} ,\;\; \alpha\in[0,1],
\end{align} 

By allowing weight $\alpha$ to vary on $\Re$, we consider general class of isotonic regression estimators
\begin{align}
	\mathcal{D}_{1,\nu}&=\{\delta_{1,c_{0,1},\alpha}^{(\nu)}:-\infty <\alpha<\infty\} \label{eq:1.9}\\\label{eq:1.10}
	\text{and}\quad \mathcal{D}_{2,\beta}&=\{\delta_{2,c_{0,2},\alpha}^{(\beta)}:-\infty <\alpha<\infty\},
\end{align}
for estimation of $\theta_1$ and $\theta_2$, respectively, when it is known apriori that $\underline{\theta}\in\Theta_0$; here $\nu$ and $\beta$ are fixed constants that are to be chosen appropriately based on the application at hand, as described above. In the literature estimators $\delta_{1,c_{0,1},\alpha}^{(\nu)}$ and $\delta_{2,c_{0,2},\alpha}^{(\beta)}$, defined by \eqref{eq:1.7} and \eqref{eq:1.8}, are also called mixed estimators of $\theta_1$ and $\theta_2$, respectively. Isotonic regression estimators have been studied by many researchers for specific probability models having independent marginals. It may be worth mentioning here that, for many probability models (e.g. independent Gaussian, independent gamma, etc.), the restricted maximum likelihood estimators (under restriction $\underline{\theta}\in \Theta_0$) are the same as isotonic regression estimators based on some LEE/SEE.

Kumar and Sharma (\citeyear{MR981031}) dealt with simultaneous estimation of ordered means of two normal distributions having known variances and characterized admissible estimators among isotonic regression estimators based on BLEEs of $\theta_1$ and $\theta_2$. For component-wise estimation of order restricted means of $k$ ($\geq 2$) normal distributions having known variances, Lee (\citeyear{MR615447}) established component-wise dominance of restricted maximum likelihood estimators over the unrestricted maximum likelihood estimators under the squared error loss function. Kelly (\citeyear{MR994278}) extended the results of Lee (\citeyear{MR615447}) under the criterion of stochastic dominance. Hwang and Peddada (\citeyear{MR1272076}) further extended the results of Lee (\citeyear{MR615447}) and Kelly (\citeyear{MR994278}) to general settings and proved certain stochastic dominance results. Kushary and Cohen (\citeyear{MR1029476}) and Misra and Singh (\citeyear{MR1366828}) studied component-wise estimation of order restricted location parameters of two exponential distributions, having known scale parameters, under the squared error loss function. Patra and Kumar (\citeyear{patra}) considered simultaneous estimation of ordered means of a bivariate normal distribution under the squared error loss function and discussed admissibility of isotonic regression estimators based on BLEEs of $\theta_1$ and $\theta_2$. Chang et. al. (\citeyear{MR3698503}) compared different isotonic regression estimators of the order restricted means of the bivariate normal distribution having a known correlation matrix under the criterion of stochastic dominance.\vspace*{2mm}

In this paper we aim to unify some of the above results proved for specific probability models by considering a general bivariate location or scale model, defined by \eqref{eq:1.1} or \eqref{eq:1.2}. The loss function considered is the squared error loss function, defined by \eqref{eq:1.3} or \eqref{eq:1.4}. For estimation of $\theta_1$ ($\theta_2$), under the restricted parameter space $\Theta_0$ and squared error loss functions \eqref{eq:1.3} or \eqref{eq:1.4}, we consider the class $\mathcal{D}_{1,\nu}$ ($\mathcal{D}_{2,\beta}$) of isotonic regression estimators defined in \eqref{eq:1.9} (\eqref{eq:1.10}), where $\nu$ $(\beta)$ are suitably chosen, as described above. We will characterize admissible estimators in $\mathcal{D}_{1,\nu}$ ($\mathcal{D}_{2,\beta}$) and obtain subclass of estimators that dominate the BLEE/BSEE $\delta_{1,c_{0,1}}(\underline{X})$ ($\delta_{2,c_{0,2}}(\underline{X})$). The rest of the paper is organized as follows. In Section 2, we introduce some useful notations, definitions, and results that are used later in the paper. Sections 3.1 and 3.2 (4.1 and 4.2), respectively, deal with estimation of the smaller and larger location (scale) parameters and under the loss functions $L_1$ and $L_2$, defined by (1.3) and (1.4), respectively. In Section 3.3 (4.3), we present a simulation study to demonstrate that an isotonic regression estimator based on BLEEs (BSEEs) may not uniformly dominate the BLEE (BSEE).

\section{\textbf{Some Useful Notations, Definitions and Results}}	\label{sec:2}
The following notations will be used throughout the paper:
\begin{itemize}
	\item[$\ast$] $\Re=(-\infty,\infty);\,\Re_{++}=(0,\infty),\,\Re^2=\Re\times \Re,\,\Re_{++}^2=\Re_{++}\times\Re_{++}$;
	\item[$\ast$] For any positive integer $m$, $\Re^m$ will denote the $m$-dimensional Euclidean space;
	\item[$\ast$] For any set $A\subseteq\Re$, $I_A(\cdot)$ will denote its indicator function;
	\item[$\ast$] $N(\mu,\sigma^2):$ Normal distribution with mean $\mu\in\Re$ and standard deviation $\sigma>0$;
	\item[$\ast$] $N_2(\mu_1,\mu_2,\sigma_1^2,\sigma_2^2,\rho):$ Bivariate normal distribution with mean vector $\begin{pmatrix} 
		\mu_1 \\ 
		\mu_2
	\end{pmatrix}\in\Re^2$ and positive definite variance-covariance matrix $\Sigma=\begin{pmatrix} 
		\sigma_1^2 & \rho \sigma_1 \sigma_2 \\ 
		\rho \sigma_1 \sigma_2  & \sigma_2^2
	\end{pmatrix} $;
	\item[$\ast$] $\Phi(z):$ distribution function of $N(0,1)$;
	\item[$\ast$] $\phi(z):$ probability density function of $N(0,1)$;
	\item[$\ast$] $=^d$ would mean equality in distribution.
	
\end{itemize}

We extend the usual orders "$\leq$" and "$<$" on $\Re$ to $\Re\cup\{-\infty,\infty\}$ with the convention that, for $a>(<)\,0,\;\frac{a}{0}=\infty\,(-\infty),\;\frac{a}{\infty}=0,\;\frac{0}{\infty}=0$ and that, for any $b\in\Re$, $b>-\infty$ and $b<\infty$.\vspace*{2mm}

Whenever, we say that an estimator $\delta_1$ dominates (or improves on) another estimator $\delta_2$ it is in the non-strict sense, i.e., it simply means that the risk of $\delta_1$ is never larger than that of $\delta_2$, for any configuration of parametric values. Also, when we say that a function $ k: A \rightarrow \Re $, where $ A \subseteq \Re $, is increasing (decreasing) it means that it is non-decreasing (non-increasing).\vspace*{2mm}

We first introduce definitions of some stochastic orders relevant to our study (see Shaked and Shanthikumar (\citeyear{MR2265633})).
\\~\\\textbf{Definition 2.1} Let $X$ and $Y$ be random variables with the Lebesgue pdfs $f(\cdot)$ and $g(\cdot)$, respectively, dfs $F(\cdot)$ and $G(\cdot)$, respectively, and distributional supports $S_X$ and $S_Y$, respectively. We say that the random variable (r.v.) $X$ is smaller than the r.v. $Y$ in the:
\\~\\(i) usual stochastic order (written as $ X \underset{\mathrm{st}}{\leq} Y$) if $F(x)\geq G(x)$, $\forall \; x\in\Re$;
\\~\\(ii)  likelihood ratio order (written as $ X \underset{\mathrm{lr}}{\leq} Y$) if 
$$ f(x) g(y)\geq f(y)g(x),\;\; \text{whenever} \; -\infty<x<y<\infty,$$
or equivalently if $S_X \subseteq S_Y$ and $\frac{g(x)}{f(x)}$ is increasing in $x\in S_Y$;
\\~\\It is well known that (see Shaked and Shanthikumar (\citeyear{MR2265633})) $ X \underset{\mathrm{lr}}{\leq} Y \Rightarrow X \underset{\mathrm{st}}{\leq} Y \iff E[h(X)]\leq E[h(Y)]$, for any increasing function $h:S_X\cup S_Y\rightarrow\Re$ for which expectations exists.\vspace*{2mm}

Now we will discuss some properties of log-concave and log-convex functions (see Pecaric et. al. (\citeyear{MR1162312})), as they are relevant to our study. We begin with the definitions of log-concave and log-convex functions.
\\~\\ \textbf{Definition 2.2} Let $D$ be a convex subset of $\Re^m$. A function $h:D\rightarrow [0,\infty)$ is said to be log-concave (log-convex) if, for all $\underline{x},\,\underline{y}\in D$ and all $\alpha\in[0,1]$,
$$h(\alpha \underline{x}+(1-\alpha) \underline{y})\geq (\leq) (h(\underline{x}))^{\alpha} (h(\underline{y}))^{1-\alpha}.$$
The following properties are well known.
\\~\\ \textbf{P1.} If $g:\Re\rightarrow [0,\infty)$ is a pdf with interior of support $(a,b)$ $(-\infty\leq a< b\leq \infty)$, then $g$ is log-concave (log-convex) on $(a,b)$ if, and only if,
\begin{equation}\label{eq:2.1}
	g(x_1)\,g(x_2-\delta)\geq (\leq) g(x_1-\delta) \,g(x_2),
\end{equation}
for all $0<\delta <b-a$ and for all $a+\delta<x_1<x_2<b$ (see Pecaric et. al. (\citeyear{MR1162312})).
\\~\\ \textbf{P2.} Let $h:\Re^m\times \Re^n\rightarrow [0,\infty)$ be a log-concave function. Then the function
$$g(\underline{x})=\int_{\Re^n} h(\underline{x},\underline{y})\,d\underline{y}$$
is log-concave on $\Re^m$ (see Prekopa (\citeyear{MR315079}) and Pecaric et. al. (\citeyear{MR1162312})).
\\~\\ \textbf{P3.} Let $-\infty\leq a<b\leq \infty$ and $-\infty\leq c<d\leq \infty$. Let $h:(a,b)\times(c,d)\rightarrow [0,\infty)$ be such that, for each $y\in (c,d)$, $h(x,y)$ is log-convex in $x\in(a,b)$. Then
$$g(x)=\int_{a}^{b} h(x,y)dy$$
is log-convex on $(a,b)$ (Artin (\citeyear{artin}) and Marshall and Olkin (\citeyear{MR552278})).\vspace*{2mm}

The following result taken from Misra and van der Meulen (\citeyear{MR1985890}) will be used in proving the main results of the paper.
\\~\\ \textbf{Proposition 2.1} Let $T_1$ and $T_2$ be random variables having distributional supports $B_1$ and $B_2$, respectively. Let $k_i:B_1\cup B_2\rightarrow \Re,\;i=1,2,$ be given functions and let $\beta_i=E[k_i(T_i)],\;i=1,2$. Suppose that $T_1\leq_{st} T_2$ ($T_2\leq_{st} T_1$). Then
\\ \textbf{(i)} $\beta_1\leq \beta_2$, provided $k_1(x)\leq k_2(x)$, $\forall\; x\in B_1\cup B_2$, and $k_1(x)$ or $k_2(x)$ is a increasing (decreasing) function of $x\in B_1\cup B_2$;
\\ \textbf{(ii)} $\beta_2\leq \beta_1$, provided $k_2(x)\leq k_1(x)$, $\forall\; x\in B_1\cup B_2$, and $k_1(x)$ or $k_2(x)$ is a decreasing (increasing) function of $x\in B_1\cup B_2$.
\vspace*{2mm}

The following result, famously known as Chebyshev's inequality, will also be used in our study (see Marshall and Olkin (\citeyear{MR2363282})).
\\~\\ \textbf{Proposition 2.2} Let $S$ be r.v. and let $k_1(\cdot)$ and $k_2(\cdot)$ be real-valued monotonic functions defined on the distributional support of the r.v. $S$. If $k_1(\cdot)$ and $k_2(\cdot)$ are monotonic functions of the same (opposite) type, then
$$E[k_1(S)k_2(S)]\geq (\leq ) E[k_1(S)] E[k_2(S)],$$
provided the above expectations exist.
\vspace*{2mm}

In the following section we consider isotonic regression estimators of order restricted location parameters $\theta_1$ and $\theta_2$ ($\theta_1\leq \theta_2$).

\section{\textbf{Isotonic Regression Estimators For Component-wise Estimation of Order Restricted Location Parameters}}	\label{sec:3}

Consider estimation of $\theta_i\;(i=1,2)$ under the location model (1.1) and the SEL function (1.3), when it is known apriori that $\underline{\theta}\in \Theta_0=\{(x_1,x_2)\in \Re^2:x_1\leq x_2\}$. The risk function of an estimator $\delta_i$ of $\theta_i$ is given by
$$R_i(\underline{\theta},\delta_i)=E_{\theta}[(\delta_i(\underline{X})-\theta_i)^2], \; \underline{\theta}\in\Theta_0,\; i=1,2.$$
Let $\delta_{i,c_{0,i}}(\underline{X})=X_i-c_{0,i},\;i=1,2,$ be the BLEE of $\theta_i,$ where $c_{0,i}$ is defined by \eqref{eq:1.5}, $i=1,2$. For suitably chosen $\nu \in\Re$ and $\beta\in\Re$, let $\mathcal{D}_{1,\nu}$ and $\mathcal{D}_{2,\beta}$ be the classes of isotonic regression estimators of $\theta_1$ and $\theta_2$ defined by (1.7)-(1.10), where $\delta_{1,\beta}(\underline{X})=X_1-\beta$ and $\delta_{2,\nu}(\underline{X})=X_2-\nu$. The choices of $\nu$ and $\beta$ are to be such that LEEs $\delta_{2,\nu}(\underline{X})=X_2-\nu$ and $\delta_{1,\beta}(\underline{X})=X_1-\beta$ are close to BLEEs $\delta_{2,c_{0,2}}(\underline{X})=X_2-c_{0,2}$ and $\delta_{1,c_{0,1}}(\underline{X})=X_1-c_{0,1}$, respectively, and dominance of $\delta_{1,c_{0,1,\alpha}}^{(\nu)}(\underline{X})$ and $\delta_{2,c_{0,2},\alpha}^{(\beta)}(\underline{X})$ over the BLEEs $\delta_{1,c_{0,1}}$ and $\delta_{2,c_{0,2}}$, respectively, is ensured for some $\alpha$'s. Our goal is to characterize admissible estimators within classes $\mathcal{D}_{1,\nu}$ and $\mathcal{D}_{2,\beta}$ of isotonic regression estimators of $\theta_1$ and $\theta_2$, respectively, and to find estimators in these classes that dominate the BLEEs $\delta_{1,c_{0,1}}(\underline{X})=X_1-c_{0,1}$ and $\delta_{2,c_{0,2}}(\underline{X})=X_2-c_{0,2}$, respectively. The following subsection deals with estimation of the location parameter $\theta_1$.

\subsection{\textbf{Isotonic Regression Estimators of the smaller Location Parameter $\theta_1$}}
\label{sec:3.1}

\noindent
\vspace*{2mm}

Let $\nu$ be a fixed real constant. Consider the class $\mathcal{D}_{1,\nu}$ of isotonic regression estimators defined by (1.7) and (1.9), where $\delta_{1,c_{0,1}}(\underline{X})=X_1-c_{0,1}$ and $\delta_{2,\nu}(\underline{x})=X_2-\nu$. Note that the BLEE $\delta_{1,c_{0,1}}(\underline{X})=\delta_{1,c_{0,1},1}^{(\nu)}(\underline{X})\in\mathcal{D}_{1,\nu}$. Define $D=X_2-X_1$, $\lambda=\theta_2-\theta_1$ ($\lambda\geq 0$), $Z_1=X_1-\theta_1$, $Z_2=X_2-\theta_2$, $Z=Z_2-Z_1$ and $\underline{Z}=(Z_1,Z_2)$. For estimating $\theta_1$ under the SEL function (1.3), the risk function of the estimator $\delta_{1,c_{0,1},\alpha}^{(\nu)}\in \mathcal{D}_{1,\nu}$,\; $-\infty<\alpha<\infty$, is given by
\begin{multline}\label{eq:3.1}
	R_1(\underline{\theta},\delta_{1,c_{0,1},\alpha}^{(\nu)})
	=E_{\underline{\theta}}[(Z_1-c_{0,1})^2 I_{(\nu-c_{0,1}-\lambda,\infty)}(Z)]\\
	+ E_{\underline{\theta}}[((Z_2+\lambda-\nu)-\alpha (Z+\lambda+c_{0,1}-\nu))^2 I_{(-\infty,\nu-c_{0,1}-\lambda)}(Z)].
\end{multline}
The risk function $R_1(\underline{\theta},\delta_{1,c_{0,1},\alpha}^{(\nu)})$ depends on $\underline{\theta}\in\Theta_0$ only through $\lambda=\theta_2-\theta_1\;(\geq 0)$. Let $f_Z(\cdot)$ denote the pdf of $Z$, so that $f_Z(z)=\int_{-\infty}^{\infty} f(s,s+z)\,ds,\; -\infty<z<\infty$. 
Clearly, for any fixed $\underline{\theta}\in\Theta_0$ (or $\lambda\geq 0)$, $R_1(\underline{\theta},\delta_{1,c_{0,1},\alpha}^{(\nu)})$ is minimized at $\alpha= \alpha_1(\lambda)$, where, for $\lambda \geq 0$,
\begin{equation}\label{eq:3.2}
	\alpha_1(\lambda)=\frac{E_{\underline{\theta}}[(Z_2+\lambda-\nu)(Z+\lambda+c_{0,1}-\nu) I_{(-\infty,\nu-c_{0,1}-\lambda)}(Z)]}{E_{\underline{\theta}}[(Z+\lambda+c_{0,1}-\nu)^2 I_{(-\infty,\nu-c_{0,1}-\lambda)}(Z)]}
	=1+\alpha_1^{*}(\lambda)
\end{equation}
\begin{equation}\label{eq:3.3}
	\text{and}\;\;\;	\alpha_1^{*}(\lambda)=\frac{E_{\underline{\theta}}[(Z_1-c_{0,1})(Z+\lambda+c_{0,1}-\nu) I_{(-\infty,\nu-c_{0,1}-\lambda)}(Z)]}{E_{\underline{\theta}}[(Z+\lambda+c_{0,1}-\nu)^2 I_{(-\infty,\nu-c_{0,1}-\lambda)}(Z)]},\;\;\lambda\geq 0.
\end{equation}
Let $S_Z$ denote the distribution support of $Z$ and
$$\psi_1(z)=E[Z_1|Z=z]=\frac{\int_{-\infty}^{\infty}s \,f(s,s+z)ds}{\int_{-\infty}^{\infty}f(s,s+z)ds},\; z\in S_Z.$$
Then, for $\lambda\geq 0$,
\begin{align*}
	\alpha_1^{*}(\lambda)
	= \frac{\int_{-\infty}^{0} (\psi_1(z+\nu-c_{0,1}-\lambda)-c_{0,1})z\,f_Z(z+\nu-c_{0,1}-\lambda) dz}{\int_{-\infty}^{0} z^2\,f_Z(z+\nu-c_{0,1}-\lambda) dz}
	=E_{\lambda}[k_1(S_{1,\lambda},\lambda)],
\end{align*}
where $k_1(z,\lambda)=\frac{\psi_1(z+\nu-c_{0,1}-\lambda)-c_{0,1}}{z},\; z<0,\; \lambda\geq 0,$ and $S_{1,\lambda}$ is a r.v. having the pdf
\begin{equation}\label{eq:3.4}
	h_{1,\lambda}(z)=\begin{cases}
		\frac{z^2 f_Z(z+\nu-c_{0,1}-\lambda)}{\int_{-\infty}^{0} s^2 f_Z(s+\nu-c_{0,1}-\lambda)\,ds},& \text{if}\;\; z<0\\
		0, &\text{otherwise}
	\end{cases},\;\;\lambda\geq 0.
\end{equation}

Using (2.1), it is easy to check if $f_Z(\cdot)$ is log-concave (log-convex) on $(-\infty,\nu-c_{0,1})$, then $S_{1,\lambda_1}\leq_{lr}S_{1,\lambda_2}$ ($S_{1,\lambda_2}\leq_{lr}S_{1,\lambda_1}$) and, consequently, $S_{1,\lambda_1}\leq_{st}S_{1,\lambda_2}$ ($S_{1,\lambda_2}\leq_{st}S_{1,\lambda_1}$), whenever $0\leq \lambda_1<\lambda_2<\infty$. \vspace*{2mm}

\noindent
The following lemma will be useful in proving the main result of this subsection.
\\~\\ \textbf{Lemma 3.1.1 (a)} Suppose that $f_Z(\cdot)$ is log-concave (log-convex) on $(-\infty,\nu-c_{0,1})$,  $\psi_1(z)=E[Z_1|Z=z]$ is increasing on $(-\infty,\nu-c_{0,1})$ and, for every $\lambda\geq 0$, $k_1(z,\lambda)=\frac{\psi_1(z+\nu-c_{0,1}-\lambda)-c_{0,1}}{z}$ is increasing (decreasing) in $z\in(-\infty,0)$. Then $\alpha_1^{*}(\lambda)$ (and hence $\alpha_1(\lambda)$), defined by (3.3) ((3.2)), is an increasing function of $\lambda\in[0,\infty)$,
$\inf_{\lambda\geq 0} \alpha_1(\lambda)= 1+\alpha_1^{*}(0)=\alpha_{1,0},$ say, and
$\sup_{\lambda\geq 0} \alpha_1(\lambda)= 1+\lim_{\lambda \to \infty}\alpha_1^{*}(\lambda)=\alpha_{1,\infty}, \text{ say}.$
\\~\\\textbf{(b)}  Suppose that $f_Z(\cdot)$ is log-concave (log-convex) on $(-\infty,\nu-c_{0,1})$,  $\psi_1(z)$ is decreasing on $(-\infty,\nu-c_{0,1})$ and, for every $\lambda\geq 0$, $k_1(z,\lambda)=\frac{\psi_1(z+\nu-c_{0,1}-\lambda)-c_{0,1}}{z}$ is decreasing (increasing) in $z\in(-\infty,0)$. Then $\alpha_1^{*}(\lambda)$ (and hence $\alpha_1(\lambda)$) is a decreasing function of $\lambda\in[0,\infty)$,
$\inf_{\lambda\geq 0} \alpha_1(\lambda) = 1+\lim_{\lambda \to \infty}\alpha_1^{*}(\lambda)=\alpha_{1,\infty}, \text{ say, and }
\sup_{\lambda\geq 0} \alpha_1(\lambda) = 1+\alpha_1^{*}(0)=\alpha_{1,0}, \text{ say}.$
\begin{proof} \textbf{(a)} Let $0\leq \lambda_1<\lambda_2 <\infty$, $T_i=S_{1,\lambda_i},\;i=1,2,$ and $m_{1,i}(z)=k_1(z,\lambda_i)=\frac{\psi_1(z+\nu-c_{0,1}-\lambda_i)-c_{0,1}}{z},\;z<0,\;i=1,2,$ where $S_{1,\lambda}$ is a r.v. having pdf given by \eqref{eq:3.4}. Then $\alpha_1^{*}(\lambda_i)=E_{\lambda_i}[m_{1,i}(T_i)],\;i=1,2$. The hypothesis that $f_Z(\cdot)$ is log-concave on $(-\infty,\nu-c_{0,1})$ implies that $T_1\leq_{lr} T_2$ ($T_2\leq_{lr} T_1$) and, consequently, $T_1\leq_{st} T_2$ ($T_2\leq_{st} T_1$). Also, under the hypothesis of (a), $m_{1,1}(t) \leq (\geq) \,m_{1,2}(t),\;\forall \; t<0$, and $m_{1,i}(t)$ is an increasing (decreasing) function of $t\in(-\infty,0)$. Now, using Proposition 2.1 (a), it follows that
	$$\alpha_1^{*}(\lambda_1)=E_{\lambda_1}[m_{1,1}(T_1)]\leq E_{\lambda_2}[m_{1,2}(T_2)]=\alpha_1^{*}(\lambda_2).$$
	\textbf{(b)} Using Proposition 2.1 (ii), the assertion follows on the lines of the proof of assertion (a) above.\vspace*{2mm}
\end{proof}

\noindent
Now we prove the main result of this subsection
\\~\\ \textbf{Theorem 3.1.1 (a)} Suppose that assumptions of Lemma 3.1.1 (a) hold. Then the estimators that are admissible within the class $\mathcal{D}_{1,\nu}$ are $\{\delta_{1,c_{0,1},\alpha}^{(\nu)}:\alpha\in[\alpha_{1,0},\alpha_{1,\infty}]\}$. Moreover, for $-\infty<\alpha_1<\alpha_2\leq \alpha_{1,0}$ or $\alpha_{1,\infty}\leq \alpha_2<\alpha_1$, the estimator $\delta_{1,c_{0,1},\alpha_2}^{(\nu)}$ dominates the estimator $\delta_{1,c_{0,1},\alpha_1}^{(\nu)}$, for any $\underline{\theta}\in\Theta_0$.
\\~\\ \textbf{(b)} Suppose that assumptions of Lemma 3.1.1 (b) hold. Then the estimators that are admissible within the class $\mathcal{D}_{1,\nu}$ are $\{\delta_{1,c_{0,1},\alpha}^{(\nu)}:\alpha\in[\alpha_{1,\infty},\alpha_{1,0}]\}$. Moreover, for $-\infty<\alpha_1<\alpha_2\leq \alpha_{1,\infty}$ or $\alpha_{1,0}\leq \alpha_2<\alpha_1<\infty$, the estimator $\delta_{1,c_{0,1},\alpha_2}^{(\nu)}$ dominates the estimator $\delta_{1,c_{0,1},\alpha_1}^{(\nu)}$, for any $\underline{\theta}\in\Theta_0$.
\begin{proof}\textbf{(a)} Let $\alpha_1(\lambda)$ be as defined by (3.2). Note that, for any fixed $\underline{\theta}\in\Theta_0$ (or fixed $\lambda\geq 0$), the risk function $R_1(\underline{\theta},\delta_{1,c_{0,1},\alpha}^{(\nu)})$, given by (3.1), is uniquely minimized at $\alpha=\alpha_1(\lambda)$, it is a strictly decreasing function of $\alpha$ on $(-\infty,\alpha_1(\lambda))$ and strictly increasing function of $\alpha$ on $(\alpha_1(\lambda),\infty)$. Since, for any $\lambda \geq 0$, $\alpha_1(\lambda)$ is a continuous function of $\lambda\in[0,\infty)$, it assumes all values between $\alpha_{1,0}=\inf_{\lambda\geq 0} \alpha_1(\lambda)$ and $\alpha_{1,\infty}=\sup_{\lambda\geq 0}\alpha_1(\lambda)$ as $\lambda$ varies on $[0,\infty)$. It follows that, each $\alpha\in[\alpha_{1,0},\alpha_{1,\infty}]$ uniquely minimizes the risk function $R_1(\underline{\theta},\delta_{1,c_{0,1},\alpha}^{(\nu)})$ at some $\underline{\theta}\in\Theta_0$ (or at some $\lambda\geq 0$). This proves that the estimators $\{\delta_{1,c_{0,1},\alpha}^{(\nu)}:\alpha\in[\alpha_{1,0},\alpha_{1,\infty}]\}$ are admissible among the estimators in the class $\mathcal{D}_{1,\nu}$. Since $\alpha_{1,0}\leq \alpha_1(\lambda)\leq \alpha_{1,\infty},\;\forall\;\lambda\geq 0$, the foregoing discussion implies that, for any $\underline{\theta}\in\Theta_0$ (or for any $\lambda\geq 0$), $R_1(\underline{\theta},\delta_{1,c_{0,1},\alpha}^{(\nu)})$ is a decreasing function of $\alpha$ on $(-\infty,\alpha_{1,0}]$ and, for any $\underline{\theta}\in\Theta_0$, it is an increasing function of $\alpha$ on $[\alpha_{1,\infty},\infty)$. This establishes the second assertion of (a).
	\\~\\ \textbf{(b)} Similar to the proof of (a), and hence omitted.
\end{proof}

In many applications verifying validity of assumptions of the above theorem will be straightforward. In this direction, the following points are noteworthy:
\begin{itemize}
	\item[$\ast$] If $f:\Re^2\rightarrow [0,\infty)$ be log-concave on $\Re^2$,
	then
	$f_Z(z)=\int_{-\infty}^{\infty} f(s,s+z)ds,\; z\in\Re,$
	is log-concave on $\Re$ (see property P2);
	\item[$\ast$] If, for every fixed $s\in\Re$, $f(s,z)$ is a log-convex function of $z\in (-\infty,\nu-c_{0,1})$, then $f_Z(z)=\int_{-\infty}^{\infty} f(s,s+z)ds,\;\;\; z\in\Re,$
	is log-convex on $(-\infty,\nu-c_{0,1})$ (see property P3);
	\item[$\ast$] Let $X_1$ and $X_2$ be independent r.v.s with pdfs $f_1(\cdot)$ and $f_2(\cdot)$, respectively, so that $f(z_1,z_2)=f_1(z_1) f_2(z_2),\; -\infty<z_1,z_2<\infty$. If $f_1$ and $f_2$ are log-concave on $\Re$, then
	$f_Z(z)=\int_{-\infty}^{\infty} f(s,s+z)ds,\, z\in\Re,$
	is log-concave on $\Re$. Similarly, if for any fixed $s\in\Re$, $f_2(s+z)$ is log-convex on (-$\infty,\nu-c_{0,1}$), then $f_Z(z)=\int_{-\infty}^{\infty} f(s,s+z)ds,\; z\in\Re,$
	is log-convex on $(-\infty,\nu-c_{0,1})$ (follows from the above two points);
	\item[$\ast$] Suppose that $\frac{f(s,s+z_2)}{f(s,s+z_1)}$ is increasing (decreasing) in $s\in\Re$, whenever $-\infty<z_1<z_2\leq \nu-c_{0,1}$. Then $\psi_1(z)=E[Z_1|Z=z]$ is an increasing (decreasing) function of $z\in (-\infty,\nu-c_{0,1})$. \vspace*{2mm}

	\noindent
	To establish the last point, note that
	$$\psi_1(z)=\frac{\int_{-\infty}^{\infty} sf(s,s+z)ds}{\int_{-\infty}^{\infty} f(s,s+z)ds}= E[Y_z],\;\; z\in(-\infty,\nu-c_{0,1}),$$
	where $Y_z,\; z\in(-\infty,\nu-c_{0,1})$, is a r.v. having the pdf
	$$h_z(s)=\frac{f(s,s+z)}{\int_{-\infty}^{\infty}f(s,s+z)ds},\; -\infty<s<\infty,\; z\in(-\infty,\nu-c_{0,1}).$$
	Let $-\infty<z_1<z_2\leq \nu-c_{0,1}$, $T_i=Y_{z_i},\;i=1,2,$ and $k_i(t)=t,\;-\infty<t<\infty,\;i=1,2$. Then the hypotheses of assertion implies that $T_1\leq_{lr}\,(\geq_{lr} T_2)$. Now, using Proposition 2.1, it follows that $\psi_1(z_1)=E[T_1]\leq E[T_2]=\psi_1(z_2)$.
\end{itemize}
\vspace*{2mm}

\noindent
Note that $\delta_{1,c_{0,1}}(\underline{X})=\delta_{1,c_{0,1},1}^{(\nu)}(\underline{x})=X_1-c_{0,1}$ is the BLEE of $\theta_1$. Using Theorem 3.1.1 (a) ((b)), one can obtain a class of estimators that dominate the BLEE $\delta_{1,c_{0,1}}^{(\nu)}$, provided $\alpha_{1,0}> (< ) \,1$. Note that $\alpha_{1,0}> (<)\,1$ if, and only if, $E_{\underline{\theta}}[(\psi_1(Z)-c_{0,1})(Z+c_{0,1}-\nu)I_{(-\infty,\nu-c_{0,1})}(Z)]> (<) \,0$. Let $k_1(s)=\psi_1(s)-c_{0,1}$ and $k_2(s)=(s+c_{0,1}-\nu) I_{(-\infty,\nu-c_{0,1})}(s)$. Then $k_2(s)$ is an increasing function of $s\in\Re$.
Using Proposition 2.2, it follows that
\begin{multline*}
	E[(\psi_1(Z)-c_{0,1})(Z+c_{0,1}-\nu) I_{(-\infty,\nu-c_{0,1})}(Z)]\\ \geq (\leq) 
	E[(\psi_1(Z)-c_{0,1})]E[(Z+c_{0,1}-\nu) I_{(-\infty,\nu-c_{0,1})}(Z)]=0,
\end{multline*}
provided $\psi_1(z)$ is an increasing (decreasing) function of $z\in (-\infty,\nu-c_{0,1})$. Thus, we have the following proposition.
\\~\\ \textbf{Proposition 3.1.1} Suppose that $\psi_1(z)=E[Z_1|Z=z]$ is increasing (decreasing) in $z\in (-\infty,\nu-c_{0,1})$. Then, $\alpha_{1,0}\geq (\leq )\; 1$.
\\~\\ The above proposition suggests that, under the assumptions of Theorem 3.1.1 (a) (or (b)), the BLEE $\delta_{1,c_{0,1}}(\underline{X})=X_1-c_{0,1}$ is inadmissible for estimating $\theta_1$ and dominating estimators are $\delta_{1,c_{0,1},\alpha_{1,0}}^{(\nu)}$, for any fixed $\nu$, provided, for any fixed $z\in S_Z$, $\psi_1(z)=E[Z_1\vert Z=z]$ is an increasing (decreasing) function of $z\in S_Z$.	\vspace*{1.5mm}

\noindent Now we illustrate some applications of Theorem 3.1.1.
\\~\\ \textbf{Example 3.1.1} Let $\underline{X}=(X_1,X_2)$ have the Lebesgue pdf belonging to location family (1.1), where $\underline{\theta}\in \Theta_0$. Assume that $f(z_1,z_2)=f(z_2,z_1)$, $\forall\; (z_1,z_2)\in\Re^2$ and $f(z_1,z_2)=f(-z_1,-z_2)\;\forall\; (z_1,z_2)\in\Re^2$. Further suppose that $f(z_1,z_2)$ is log-concave on $\Re^2$. Then, by property P2, $f_Z(z)$ is log-concave on $\Re$. Consider estimation of $\theta_1$ under the SEL function $L_1$, defined by (1.3), when it is known apriori that $\underline{\theta}\in\Theta_0$. Then $Z_i=^{d} -Z_i,\;i=1,2,$ and $(Z_1,Z_2)=^{d} (-Z_2,-Z_1)$. Consequently $c_{0,1}=c_{0,2}=0$ and, for $-\infty<z<\infty$, 
\begin{align*}
	\psi_1(z)\!=\!E_{\underline{\theta}}[Z_1|Z_2\!-\!Z_1\!=\!z]\!
	=\!-z\!+\!E_{\underline{\theta}}[Z_2|Z_2\!-\!Z_1\!=\!z]
	\!=\!-z\!-\!\psi_1(z)\;\; ((Z_1,Z_2)=^{d} (-Z_2,-Z_1)),
\end{align*}
i.e., $ \psi_1(z)=-\frac{z}{2},\;\; -\infty<z<\infty$. Thus, $\psi_1(z)$ is a decreasing function of $z\in\Re$. Also, for any $\lambda\geq 0$,
$$k_1(z,\lambda)=\frac{\psi_1(z+\nu-c_{0,1}-\lambda)-c_{0,1}}{z}=-\frac{1}{2} +\frac{\lambda-\nu}{2z}$$
is a decreasing function of $z\in(-\infty,0)$ for the choice $\nu\leq 0$. As $c_{0,2}=0$, here, an appropriate choice of $\nu$ is 0, which amounts to considering isotonic regression of $(X_1,X_2)$. For $\nu=0$, we have $\alpha_{1,0}=\alpha_1(0)=\frac{1}{2}$. Using Theorem 3.1.1 (b), it follows that the estimators $\{\delta_{1,0,\alpha}^{(0)}:\alpha_{1,\infty}\leq \alpha\leq \frac{1}{2}\}$ are admissible among the isotonic regression estimators in the class $\mathcal{D}_{1,0}=\{\delta_{1,0,\alpha}^{(0)}:-\infty< \alpha<\infty\}$, where
$$ \delta_{1,0,\alpha}^{(0)}(\underline{X})=\begin{cases}
	X_1, & \text{ if} \;\;\; X_1\leq X_2\\
	\alpha X_1+(1-\alpha)X_2, & \text{ if}\;\;\; X_1>X_2
\end{cases}, \;\; -\infty<\alpha<\infty.$$
Also, for any $\alpha_1>\alpha_2 \geq \frac{1}{2}$, the estimator $\delta_{1,0,\alpha_2}^{(0)}(\underline{X})$ dominates the estimator $\delta_{1,0,\alpha_1}^{(0)}(\underline{X})$. In particular the estimator
$$\delta_{1,0,0.5}^{(0)}(\underline{X})=\begin{cases}
	X_1,&\text{ if} \;\;X_1\leq X_2\\
	\frac{X_1+X_2}{2}, &\text{ if} \;\; X_1>X_2
\end{cases}=\min\bigg\{X_1,\frac{X_1+X_2}{2}\bigg\}$$
dominates the BLEE $\delta_{1,0}(\underline{X})=\delta_{1,0,1}^{(\nu)}=X_1$.\vspace*{2mm}

Examples of log-concave densities $f(z_1,z_2)$ satisfying the assumptions of this example are:
\begin{itemize}
	\item $ f(z_1,z_2) =\frac{1}{2 \pi \sigma^2 \sqrt{1-\rho^2}} e^{-\frac{1}{2(1-\rho^2)\sigma^2}\left[z_1^2-2 \rho \, z_1 z_2+z_2^2\right]}, \; (z_1,z_2)\in \Re^2,$ where $\sigma>0$ and $-1<\rho <1$.
	\\~\\ \item $f(z_1,z_2)=\frac{1}{\sigma^2}\frac{e^{-\frac{z_1}{\sigma}}}{\left(1+e^{-\frac{z_1}{\sigma}}\right)^2}\frac{e^{-\frac{z_2}{\sigma}}}{\left(1+e^{-\frac{z_2}{\sigma}}\right)^2},\;(z_1,z_2)\in \Re^2,$ where $\sigma>0$.
	\\~\\\item  $f(z_1,z_2)=\frac{1}{4\sigma^2} e^{-\frac{|z_1|}{\sigma}} e^{-\frac{|z_2|}{\sigma}} ,\;(z_1,z_2)\in \Re^2,$ where $\sigma>0$.
\end{itemize}\vspace*{2mm}

\noindent
\textbf{Example 3.1.2} Let $\underline{X}=(X_1,X_2)\sim N_2(\theta_1,\theta_2,\sigma_1^2,\sigma_2^2,\rho)$, where $-\infty<\theta_1\leq \theta_2<\infty$, $\theta_1$ and $\theta_2$ are unknown, and $\sigma_i>0,\;i=1,2,$ and $\rho\in(-1,1)$ are known. Consider estimation of $\theta_1$ under the squared error loss $L_1$, given by (1.3). Here $Z_i\sim N(0,\sigma_i^2),\;i=1,2,\; c_{0,i}=0,\;i=1,2,$ $Z_1|Z=z\sim N\left(\frac{\sigma_1(\rho\sigma_2-\sigma_1)z}{\tau^2},\frac{(1-\rho^2)\sigma_1^2 \sigma_2^2}{\tau^2}\right),\; -\infty<z<\infty,$ $Z\sim N(0,\tau^2)$ and $\psi_1(z)=E[Z_1|Z=z]=\frac{\sigma_1(\rho\sigma_2-\sigma_1)z}{\tau^2},\;-\infty<z<\infty,$ where $\tau^2=\sigma_1^2+\sigma_2^2-2\rho \sigma_1\sigma_2$. Moreover, $f_Z(z)$ is log-concave on $\Re$. Clearly $\psi_1(z)$ is decreasing (increasing) in $z\in(-\infty,0)$ if $\rho \sigma_2<\sigma_1$ ($\rho \sigma_2>\sigma_1$). Moreover, for any fixed $\lambda\geq 0$, $k_1(z,\lambda)=\frac{\sigma_1(\rho\sigma_2-\sigma_1)}{\tau^2}\frac{z+\nu-\lambda}{z}$ is decreasing (increasing) in $z\in (-\infty,0)$ provided $\rho \sigma_2<\sigma_1$ ($\rho \sigma_2>\sigma_1$) and $\nu\leq 0$. Since $c_{0,2}=0$, an appropriate choice of $\nu$ in this case is $\nu=0$. For $\nu=0$.
$$\alpha_{1,0}\!=\!1-\alpha_1^{*}(0)\!=\!1+\frac{E[\psi_1(Z)Z\,I_{(-\infty,0)}(Z)]}{E[Z^2\,I_{(-\infty,0)}(Z)]}\!=\!1+\frac{\sigma_1(\rho\sigma_2-\sigma_1)}{\tau^2}\!=\!\frac{\sigma_2(\sigma_2\!-\!\rho\sigma_1)}{\tau^2}\!=\!\beta_0,\; (\text{say})$$
and
$$ \alpha_{1,\infty}\!=\!\lim_{\lambda\to \infty}\!(1+\alpha_1^{*}(\lambda))\!=\!1+\frac{\sigma_1(\rho\sigma_2\!-\!\sigma_1)}{\tau^2} \lim_{\lambda\to \infty} \!\frac{\int_{-\infty}^{0} (z-\lambda) z\frac{1}{\tau}\phi\left(\frac{z-\lambda}{\tau}\right)dz}{\int_{-\infty}^{0}  z^2 \frac{1}{\tau}\phi\left(\frac{z-\lambda}{\tau}\right)dz}
\!=\!\begin{cases}  -\infty,\!\!\! &\text{if} \; \rho \sigma_2<\sigma_1\\
	1, \!\!\!& \text{if} \; \rho \sigma_2=\sigma_1\\
	\infty,\!\!\! & \text{if} \;\rho \sigma_2 >\sigma_1 \end{cases}.$$
Here $\delta_{1,0,\beta_0}^{(0)}(\underline{X})$
is the restricted MLE (RMLE) of $\theta_1$ (see Patra and Kumar (\citeyear{patra})). Consider the following cases.
\\~\\\textbf{Case I:} $\frac{\rho \sigma_2}{\sigma_1}<1$ \vspace*{2mm}

In this case $\alpha_{1,0}=\beta_0=\frac{\sigma_2(\sigma_2-\rho\sigma_1)}{\tau^2}<1,\; \alpha_{1,\infty}=-\infty,\;\psi_1(z)=\frac{\sigma_1(\rho\sigma_2-\sigma_1)}{\tau^2}z$ is decreasing in $z\in (-\infty,0)$ and, for every fixed $\lambda\geq 0$, $k_1(z,\lambda)=\frac{\psi_1(z-\lambda)}{z}=\frac{\sigma_1(\rho\sigma_2-\sigma_1)}{\tau^2}\frac{z-\lambda}{z}$ is decreasing in $z\in (-\infty,0)$. Using Theorem 3.1.1 (b), it follows that the estimators $\{\delta_{1,0,\alpha}^{(0)}:\alpha\in (-\infty,\beta_0]\}$ are admissible within the class $\mathcal{D}_{1,0}=\{\delta_{1,0,\alpha}^{(0)}:-\infty<\alpha<\infty\}$ of isotonic regression estimators of $\theta_1$. Moreover the estimators $\{\delta_{1,0,\alpha}^{(0)}:\alpha>\beta_0\}$ are inadmissible and, for $\beta_0\leq \alpha_2<\alpha_1$, the estimator $\delta_{1,0,\alpha_2}^{(0)}$ dominates the estimator $\delta_{1,0,\alpha_1}^{(0)}$. In particular the RMLE $\delta_{1,0,\beta_0}^{(0)}$ is admissible  within the class $\mathcal{D}_{1,0}$ of isotonic regression estimators of $\theta_1$ and it dominates the BLEE $\delta_{1,0}(\underline{X})=\delta_{1,0,1}^{(0)}=X_1$. For the special case $\rho=0$, the dominance of RMLE over the BLEE is also shown in Kubokawa and Saleh (\citeyear{MR1370413}).
\\~\\\textbf{Case II:} $\frac{\rho \sigma_2}{\sigma_1}>1$ \vspace*{2mm}

In this case $\alpha_{1,0}=\beta_0=\frac{\sigma_2(\sigma_2-\rho\sigma_1)}{\tau^2}>1,\; \alpha_{1,\infty}=\infty,\;\psi_1(z)=\frac{\sigma_1(\rho\sigma_2-\sigma_1)}{\tau^2}z$ is increasing in $z\in (-\infty,0)$ and, for every fixed $\lambda\geq 0$, $k_1(z,\lambda)=\frac{\psi_1(z-\lambda)}{z}=\frac{\sigma_1(\rho\sigma_2-\sigma_1)}{\tau^2}\frac{z-\lambda}{z}$ is increasing in $z\in (-\infty,0)$. Using Theorem 3.1.1 (a), we conclude that the estimators $\{\delta_{1,0,\alpha}^{(0)}:\alpha\in [\beta_0,\infty)\}$ are admissible within the class $\mathcal{D}_{1,0}=\{\delta_{1,0,\alpha}^{(0)}:-\infty<\alpha<\infty\}$ of isotonic regression estimators of $\theta_1$. Moreover, the estimators $\{\delta_{1,0,\alpha}^{(0)}:\alpha<\beta_0\}$ are inadmissible for estimating $\theta_1$ and, for $\infty< \alpha_1<\alpha_2\leq \beta_0$, the estimator $\delta_{1,0,\alpha_2}^{(0)}$ dominates the estimator $\delta_{1,0,\alpha_1}^{(0)}$. In particular the RMLE $\delta_{1,0,\beta_0}^{(0)}$ is admissible  within the class $\mathcal{D}_{1,0}$ of isotonic regression estimators of $\theta_1$ and it dominates the BLEE $\delta_{1,0}(\underline{X})=\delta_{1,0,1}^{(0)}=X_1$.	
\\~\\\textbf{Case III:} $\frac{\rho \sigma_2}{\sigma_1}=1$ \vspace*{2mm}

In this case $\alpha_{1,0}=\beta_0=\alpha_{1,\infty}=1$ and the BLEE $\delta_{1,0}(\underline{X})=\delta_{1,0,1}^{(0)}(\underline{X})=X_1$ is also the RMLE. Also, $\psi_1(z)=0=k(z,\lambda),\; z\in(-\infty,0),\;\lambda\geq 0$. Using Theorem 3.1.1, it follows that the BLEE/RMLE $\delta_{1,0}(\underline{X})=\delta_{1,0,1}^{(0)}(\underline{X})=X_1$ is the only admissible estimator within the class $\mathcal{D}_{1,0}$ of isotonic regression estimators of $\theta_1$. Any other isotonic regression estimator in $\mathcal{D}_{1,0}$ is dominated by the BLEE/RMLE $\delta_{1,0}(\underline{X})=\delta_{1,0,1}^{(0)}(\underline{X})=X_1$.
\\~\\	\textbf{Example 3.1.3.} Let $f(z_1,z_2)=f_1(z_1) f_2(z_2),\, \underline{z}=(z_1,z_2)\in\Re^2,$
where, for known positive constants $\sigma_1>0$ and $\sigma_2>0$,
$ f_i(z)= \frac{1}{\sigma_i} e^{-\frac{z}{\sigma_i}},\text{ if}\;z>0, i=1,2.$
Here $c_{0,i}=\sigma_i,\;i=1,2$, and the BLEE of $\theta_1$ is $\delta_{1,\sigma_1}(\underline{X})=\delta_{1,\sigma_1,1}^{(\nu)}(\underline{X})=X_1-\sigma_1$. Moreover,
$$f_Z(z)=\begin{cases}  \frac{1}{\sigma_1+\sigma_2} e^{\frac{z}{\sigma_1}},&\text{ if} \;\; z<0  \\
	\frac{1}{\sigma_1+\sigma_2} e^{-\frac{z}{\sigma_2}},&\text{ if} \;\; z\geq 0   \end{cases}$$		
is log-concave on $\Re$. The conditional pdf of $Z_1$ given $Z=z \; (-\infty<z<\infty)$ is
$f_{Z_1|Z}(z_1|z) =  \frac{\sigma_1+\sigma_2}{\sigma_1 \sigma_2}\, e^{-\left(\frac{1}{\sigma_1}+\frac{1}{\sigma_2}\right)(z_1-\max\{-z,0\})}, \text{ if} \; z_1\geq \max\{-z,0\},  $
and 
$$\psi_1(z)=E[Z_1|Z=z]=\max\{-z,0\}+\frac{\sigma_1 \sigma_2}{\sigma_1+\sigma_2} = \begin{cases}
	-z+\frac{\sigma_1 \sigma_2}{\sigma_1+\sigma_2}, &\text{ if} \;\; z<0\\
	\frac{\sigma_1 \sigma_2}{\sigma_1+\sigma_2}, & \text{ if} \;\; z\geq 0
\end{cases}$$
is decreasing in $z\in(-\infty,\nu-\sigma_1)$. For any fixed $\lambda \geq 0$,
$$k_1(z,\lambda)=\frac{\psi_1(z+\nu-\sigma_1-\lambda)-\sigma_1}{z}=\begin{cases}
	\left(\lambda+ \frac{\sigma_1\sigma_2}{\sigma_1+ \sigma_2}-\nu\right)\frac{1}{z}-1,& \text{ if} \;\; z<\lambda+\sigma_1-\nu\\
	\frac{-\sigma_1^2}{(\sigma_1+\sigma_2)z}, &\text{ if} \;\; z\geq \lambda+\sigma_1-\nu
\end{cases}$$
is decreasing in $z\in(-\infty,0)$, provided $\nu\leq \frac{\sigma_1 \sigma_2}{\sigma_1+\sigma_2}$. Since $c_{0,2}=\sigma_2$, an appropriate choice of $\nu$ in this case is $\nu=\frac{\sigma_1 \sigma_2}{\sigma_1+\sigma_2}$. For $\nu=\frac{\sigma_1 \sigma_2}{\sigma_1+\sigma_2}$, it is easy to verify that $\alpha_{1,0}=\alpha_1(0)=0$ and $\alpha_{1,\infty}=\lim_{\lambda \to \infty}\alpha_1(\lambda)=-\infty$.
Using Theorem 3.1.1 (b), we conclude that the estimators $\{\delta_{1,\sigma_1,\alpha}^{(\nu)}: \alpha\in(-\infty,0] \}$ are admissible within the class $\{\delta_{1,\sigma_1,\alpha}^{(\nu)}: -\infty<\alpha<\infty \}$ of isotonic regression estimators of $\theta_1$; here $\nu=\frac{\sigma_1 \sigma_2}{\sigma_1+\sigma_2}$. Moreover, for $0\leq \alpha_2<\alpha_1<\infty$, the isotonic regression estimator $ \delta_{1,\sigma_1,\alpha_2}^{(\nu)}$ dominates the estimator $\delta_{1,\sigma_1,\alpha_1}^{(\nu)}$. In particular the BLEE $\delta_{1,\sigma_1}(\underline{X})=\delta_{1,\sigma_1,1}^{(\nu)}(\underline{X})=X_1-\sigma_1$ is inadmissible for estimating $\theta_1$ and is dominated by the isotonic regression estimator
$$\delta_{1,\sigma_1,0}^{(\nu)}(\underline{X})=\begin{cases}
	X_1-\sigma_1,&\text{ if} \;\; X_1-\sigma_1\leq X_2-\frac{\sigma_1 \sigma_2}{\sigma_1+\sigma_2}\\
	X_2-\frac{\sigma_1 \sigma_2}{\sigma_1+\sigma_2}, &\text{ if} \; \; X_1-\sigma_1> X_2-\frac{\sigma_1 \sigma_2}{\sigma_1+\sigma_2}
\end{cases}.$$
The above estimator is a member of a class of dominating estimators proposed in Pal and Kushary (\citeyear{MR1165709}). Note that, here, isotonic regression estimators based on BLEE $(X_1-\sigma_1,X_2-\sigma_2)$ may not be able to provide improvements over the BLEE $\delta_{1,\sigma_1}(\underline{X})=X_1-\sigma_1$, rather an isotonic regression estimator based on $(X_1-\sigma_1,X_2-\frac{\sigma_1\sigma_2}{\sigma_1+\sigma_2})$ provides an improvement over the BLEE $\delta_{1,\sigma_1}(\underline{X})$ (see Section (3.3)).

\subsection{\textbf{Isotonic Regression Estimators of the Larger Location Parameter $\theta_2$}}
\label{sec:3.2}
\noindent
\vspace*{2mm}

In this section, we deal with estimation of the larger location parameter $\theta_2$ under the location model \eqref{eq:1.1} and the SEL function $L_2$, defined by \eqref{eq:1.3}, when it is known apriori that $\underline{\theta}\in\Theta_0$. We continue to follow the notations of Section 3.1. Let $\beta$ be a fixed real constant, to be suitably chosen as described in Section 1. Consider the class $\mathcal{D}_{2,\beta}$ of isotonic regression estimators of $\theta_2$, defined by (1.8) and (1.10); here $\delta_{2,c_{0,2}}(\underline{X})=X_2-c_{0,2}$ and $\delta_{1,\beta}(\underline{X})=X_1-\beta$. Note that the BLEE $\delta_{2,c_{0,2}}(\underline{X})=\delta_{2,c_{0,2},0}^{(\nu)}(\underline{X})=X_2-c_{0,2}$ is a member of the class $\mathcal{D}_{2,\beta}$ for any fixed $\beta$. Under the notations of Section 3.1, for any fixed $\underline{\theta}\in \Theta_0$ (or $\lambda\geq 0)$,
\begin{multline}\label{eq:3.5}
	R_2(\underline{\theta},\delta_{2,c_{0,2},\alpha}^{(\beta)})
	=E_{\underline{\theta}}[(Z_2-c_{0,2})^2 I_{(c_{0,2}-\beta-\lambda,\infty)}(Z)] \\
	\;\; \quad+ E_{\underline{\theta}}[((Z_2-c_{0,2})-\alpha (Z+\lambda+\beta-c_{0,2}))^2 I_{(-\infty,c_{0,2}-\beta-\lambda)}(Z)].
\end{multline}
is minimized at $\alpha= \alpha_2(\lambda)$, where
\begin{equation}\label{eq:3.6}
	\alpha_2(\lambda)=\frac{E_{\underline{\theta}}[(Z_2-c_{0,2})(Z+\lambda+\beta-c_{0,2}) I_{(-\infty,c_{0,2}-\beta-\lambda)}(Z)]}{E_{\underline{\theta}}[(Z+\lambda+\beta-c_{0,2})^2 I_{(-\infty,c_{0,2}-\beta-\lambda)}(Z)]},\;\;\lambda\geq 0.
\end{equation}
Let $\psi_2(z)=E[Z_2|Z=z]=\frac{\int_{-\infty}^{\infty}s \,f(s-z,s)ds}{\int_{-\infty}^{\infty}f(s-z,s)ds},\; z\in S_Z,$ so that
$\psi_2(z)=z+\psi_1(z),\;z\in S_Z$, where $\psi_1(z)=E[Z_1|Z=z],\;z\in S_Z$. Then, for $\lambda\geq 0$,
\begin{align}\label{eq:3.7}
	\alpha_2(\lambda)
	\!= \! \frac{\int_{-\infty}^{0} (\psi_2(z+c_{0,2}-\beta-\lambda)-c_{0,2})z\,f_Z(z+c_{0,2}-\beta-\lambda) dz}{\int_{-\infty}^{0} z^2\,f_Z(z+c_{0,2}-\beta-\lambda) dz}\!
	=\!E_{\lambda}[k_2(S_{2,\lambda},\lambda)],
\end{align}
where $k_2(z,\lambda)=\frac{\psi_2(z+c_{0,2}-\beta-\lambda)-c_{0,2}}{z},\; z<0,\; \lambda\geq 0,$ and $S_{2,\lambda}$ is a r.v. having the pdf
\begin{equation}\label{eq:3.8}
	h_{2,\lambda}(z)=\begin{cases}
		\frac{z^2 f_Z(z+c_{0,2}-\beta-\lambda)}{\int_{-\infty}^{0} s^2 f_Z(s+c_{0,2}-\beta-\lambda)\,ds},& \text{if}\;\; z<0\\
		0, &\text{otherwise}
	\end{cases},\;\;\lambda\geq 0.
\end{equation}

Using property P1, it is easy to check if $f_Z(\cdot)$ is log-concave (log-convex) on $(-\infty,c_{0,2}-\beta)$, then $S_{2,\lambda_1}\leq_{lr}S_{2,\lambda_2}$ ($S_{2,\lambda_2}\leq_{lr}S_{2,\lambda_1}$) and, consequently, $S_{2,\lambda_1}\leq_{st}S_{2,\lambda_2}$ ($S_{2,\lambda_2}\leq_{st}S_{2,\lambda_1}$), whenever $0\leq \lambda_1<\lambda_2<\infty$. \vspace*{2mm}

\noindent
The following lemma, whose proof is similar to that of Lemma 3.1.1 (a), will be instrumental in proving the main result of this subsection.
\\~\\ \textbf{Lemma 3.2.1 (a)} Suppose that $f_Z(\cdot)$ is log-concave (log-convex) on $(-\infty,c_{0,2}-\beta)$,  $\psi_2(z)=E[Z_2|Z=z]$ is increasing on $(-\infty,c_{0,2}-\beta)$ and, for every $\lambda\geq 0$, $k_2(z,\lambda)=\frac{\psi_2(z+c_{0,2}-\beta-\lambda)-c_{0,2}}{z}$ is increasing (decreasing) in $z\in(-\infty,0)$. Then $\alpha_2(\lambda)$, defined by \eqref{eq:3.7}, is an increasing function of $\lambda\in[0,\infty)$,
$\inf_{\lambda\geq 0} \alpha_2(\lambda) = \alpha_2(0)=\alpha_{2,0},  \text{ say, and }
\sup_{\lambda\geq 0} \alpha_2(\lambda) = \lim_{\lambda \to \infty}\alpha_2(\lambda)=\alpha_{2,\infty}, \text{ say}.$
\\~\\\textbf{(b)}  Suppose that $f_Z(\cdot)$ is log-concave (log-convex) on $(-\infty,c_{0,2}-\beta)$,  $\psi_2(z)$ is decreasing on $(-\infty,c_{0,2}-\beta)$ and, for every $\lambda\geq 0$, $k_2(z,\lambda)=\frac{\psi_2(z+c_{0,2}-\beta-\lambda)-c_{0,2}}{z}$ is decreasing (increasing) in $z\in(-\infty,0)$. Then $\alpha_2(\lambda)$ is a decreasing function of $\lambda\in[0,\infty)$,
$\inf_{\lambda\geq 0} \alpha_2(\lambda) = \lim_{\lambda \to \infty}\alpha_2(\lambda)=\alpha_{2,\infty}, \text{ say, and }
\sup_{\lambda\geq 0} \alpha_2(\lambda) = \alpha_2(0)=\alpha_{2,0}, \text{ say}.$
\vspace*{2mm}

\noindent
Now we have the main result of this subsection. The proof of the theorem, being similar to the proof of Theorem 3.1.1, is omitted.
\\~\\ \textbf{Theorem 3.2.1 (a)} Suppose that assumptions of Lemma 3.2.1 (a) hold. Then, the estimators that are admissible in the class $\mathcal{D}_{2,\beta}$ are $\{\delta_{2,c_{0,2},\alpha}^{(\beta)}:\alpha\in[\alpha_{2,0},\alpha_{2,\infty}]\}$. Moreover, for $-\infty<\alpha_1<\alpha_2\leq \alpha_{2,0}$ or $\alpha_{2,\infty}\leq \alpha_2<\alpha_1$, the estimator $\delta_{2,c_{0,2},\alpha_2}^{(\beta)}$ dominates the estimator $\delta_{2,c_{0,2},\alpha_1}^{(\beta)}$, for any $\underline{\theta}\in\Theta_0$.
\\~\\ \textbf{(b)} Suppose that assumptions of Lemma 3.2.1 (b) hold. Then, the estimators that are admissible in the class $\mathcal{D}_{2,\beta}$ are $\{\delta_{2,c_{0,2},\alpha}^{(\beta)}:\alpha\in[\alpha_{2,\infty},\alpha_{2,0}]\}$. Moreover, for $-\infty<\alpha_1<\alpha_2\leq \alpha_{2,\infty}$ or $\alpha_{2,0}\leq \alpha_2<\alpha_1<\infty$, the estimator $\delta_{2,c_{0,2},\alpha_2}^{(\beta)}$ dominates the estimator $\delta_{2,c_{0,2},\alpha_1}^{(\beta)}$, for any $\underline{\theta}\in\Theta_0$.
\vspace*{2mm}

Note that $\delta_{2,c_{0,2}}(\underline{X})=\delta_{2,c_{0,2},0}^{(\beta)}(\underline{X})=X_2-c_{0,2}$ is the BLEE of $\theta_2$. Using Theorem 3.2.1 (a) ((b)), one can obtain a class of estimators that dominate the BLEE $\delta_{2,c_{0,2}}$, provided $\alpha_{2,0}> (< )\, 0$. Note that $\alpha_{2,0}> (<)\,0$ if, and only if, $E_{\underline{\theta}}[(\psi_2(Z)-c_{0,2})(Z+\beta-c_{0,2})I_{(-\infty,c_{0,2}-\beta)}(Z)]> (<) \,0$. Let $k_1(s)=\psi_2(s)-c_{0,2}$ and $k_2(s)=(s+\beta-c_{0,2}) I_{(-\infty,c_{0,2}-\beta)}(s)$. Then $k_2(s)$ is an increasing function of $s\in\Re$.
Using Proposition 2.2, it follows that
\begin{multline*}
	E[(\psi_2(Z)-c_{0,2})(Z+\beta-c_{0,2}) I_{(-\infty,c_{0,2}-\beta)}(Z)]\\ \geq (\leq) 
	E[(\psi_2(Z)-c_{0,2})]E[(Z+\beta-c_{0,2}) I_{(-\infty,c_{0,2}-\beta)}(Z)]=0,
\end{multline*}
provided $\psi_2(z)$ is an increasing (decreasing) function of $z\in (-\infty,c_{0,2}-\beta)$. 
Thus, we have the following proposition.
\\~\\ \textbf{Proposition 3.2.2} Suppose that $\psi_2(z)=E[Z_2|Z=z]$ is increasing (decreasing) in $z\in (-\infty,c_{0,2}-\beta)$. Then, $\alpha_{2,0}\geq (\leq )\; 0$. 
\\~\\Now we illustrate some applications of Theorem 3.2.1.
\\~\\ \textbf{Example 3.2.1} Under the set-up of Example 3.1.1, consider estimation of $\theta_2$ under the SEL function $L_2$, defined by \eqref{eq:1.3}, when it is known apriori that $\underline{\theta}\in\Theta_0$. Then $c_{0,1}=c_{0,2}=0$ and $\psi_2(z)=z+\psi_1(z)=\frac{z}{2},\;\; -\infty<z<\infty$. Thus, $\psi_2(z)$ is an increasing function of $z\in\Re$. Also, for any $\lambda\geq 0$,
$$k_2(z,\lambda)=\frac{\psi_2(z+c_{0,2}-\beta-\lambda)-c_{0,2}}{z}=\frac{1}{2} -\frac{\beta+\lambda}{2z}$$
is an increasing function of $z\in(-\infty,0)$, for any $\beta\geq 0$. As $c_{0,1}=0$, here the suitable choice of $\beta$ is 0. That amounts to considering isotonic regression of $(X_1,X_2)$. One can verify that, for $\beta=0$,
$\alpha_{2,0}=\alpha_2(0)=\frac{1}{2}\; \text{and} \; \alpha_{2,\infty}=\lim_{\lambda\to \infty} \alpha_{2,\infty}.$ \vspace*{2mm}

Now using Theorem 3.2.1 (a), it follows that the estimators $\{\delta_{2,0,\alpha}^{(0)}:\frac{1}{2}\leq \alpha\leq\alpha_{2,\infty} \}$ are admissible among the isotonic regression estimators in the class $\mathcal{D}_{2,0}=\{\delta_{2,0,\alpha}^{(0)}:-\infty< \alpha<\infty\}$, where
$$ \delta_{2,0,\alpha}^{(0)}(\underline{X})=\begin{cases}
	X_2, & \text{ if} \;\;\; X_1\leq X_2\\
	\alpha X_1+(1-\alpha)X_2, & \text{ if}\;\;\; X_1>X_2
\end{cases}, \;\; -\infty<\alpha<\infty.$$
Also, for $-\infty<\alpha_1<\alpha_2 \leq \frac{1}{2}$, the estimator $\delta_{2,0,\alpha_2}^{(0)}(\underline{X})$ dominates the estimator $\delta_{2,0,\alpha_1}^{(0)}(\underline{X})$. In particular the estimator
$$\delta_{2,0,0.5}^{(0)}(\underline{X})=\begin{cases}
	X_2,&\text{ if} \;\;X_1\leq X_2\\
	\frac{X_1+X_2}{2}, &\text{ if} \;\; X_1>X_2
\end{cases}=\max\bigg\{X_2,\frac{X_1+X_2}{2}\bigg\}$$
dominates the BLEE $\delta_{2,0}(\underline{X})=\delta_{2,0,0}^{(0)}=X_2$.
\\~\\\textbf{Example 3.2.2.} Let $\underline{X}=(X_1,X_2)$ follow the bivariate normal distribution described in Example 3.1.2. Consider estimation of $\theta_2$ under the squared error loss $L_2$, given by (1.3), when it is known that $\underline{\theta}\in \Theta_0$. Here, for any fixed $\lambda\geq 0$, $k_2(z,\lambda)=\frac{\sigma_2(\sigma_2-\rho\sigma_1)}{\tau^2}\frac{z-\beta-\lambda}{z}$ is decreasing (increasing) in $z\in (-\infty,0)$, provided $\rho \sigma_1>\sigma_2$ ($\rho \sigma_1<\sigma_2$) and $\beta\geq 0$. Since $c_{0,1}=0$, an appropriate choice of $\beta$ in this case is $\beta=0$. For $\beta=0$, we have
$$\alpha_{2,0}\!=\!\alpha_2(0)\!=\!\frac{E[\psi_2(Z)Z\,I_{(-\infty,0)}(Z)]}{E[Z^2\,I_{(-\infty,0)}(Z)]}\!=\!1\!+\!\frac{E[\psi_1(Z)Z\,I_{(-\infty,0)}(Z)]}{E[Z^2\,I_{(-\infty,0)}(Z)]}\!=\!\frac{\sigma_2(\sigma_2-\rho\sigma_1)}{\tau^2}\!=\!\beta_0$$
and
$$ \alpha_{2,\infty}\!=\!\lim_{\lambda\to \infty}\alpha_2(\lambda)\!=\!\frac{\sigma_2(\sigma_2-\rho\sigma_1)}{\tau^2} \lim_{\lambda\to \infty} \frac{\int_{-\infty}^{0} (z-\lambda) z\frac{1}{\tau}\phi\left(\frac{z-\lambda}{\tau}\right)dz}{\int_{-\infty}^{0}  z^2 \frac{1}{\tau}\phi\left(\frac{z-\lambda}{\tau}\right)dz}
\!=\!\begin{cases}  -\infty,\!\!\! &\text{if} \, \rho \sigma_1>\sigma_2\\
	0,\!\!\! &\text{if} \, \rho \sigma_1=\sigma_2\\
	\infty,\!\!\! &\text{if} \, \rho \sigma_1 <\sigma_2 \end{cases}.$$
Here $\delta_{2,0,\beta_0}^{(0)}(\underline{X})$
is the restricted MLE (RMLE) of $\theta_2$. Consider the following cases.
\\~\\\textbf{Case I:} $\rho \sigma_1>\sigma_2$ \vspace*{2mm}

In this case $\alpha_{2,0}=\beta_0<0,\; \alpha_{2,\infty}=-\infty,\;\psi_2(z)=\frac{\sigma_2(\sigma_2-\rho\sigma_1)}{\tau^2}z$ is decreasing in $z\in (-\infty,0)$ and, for every fixed $\lambda\geq 0$, $k_2(z,\lambda)=\frac{\psi_2(z-\lambda)}{z}=\frac{\sigma_2(\sigma_2-\rho\sigma_1)}{\tau^2}\frac{z-\lambda}{z}$ is decreasing in $z\in (-\infty,0)$. Using Theorem 3.2.1 (b), it follows that the estimators $\{\delta_{2,0,\alpha}^{(0)}:\alpha\in (-\infty,\beta_0]\}$ are admissible within the class $\mathcal{D}_{2,0}=\{\delta_{2,0,\alpha}^{(0)}:-\infty<\alpha<\infty\}$ of isotonic regression estimators of $\theta_2$. Moreover, the estimators $\{\delta_{2,0,\alpha}^{(0)}:\alpha>\beta_0\}$ are inadmissible and, for $\beta_0\leq \alpha_2<\alpha_1$, the estimator $\delta_{2,0,\alpha_2}^{(0)}$ dominates the estimator $\delta_{2,0,\alpha_1}^{(0)}$. In particular the RMLE $\delta_{2,0,\beta_0}^{(0)}$ is admissible  within the class $\mathcal{D}_{2,0}$ of isotonic regression estimators of $\theta_2$ and it dominates the BLEE $\delta_{2,0}(\underline{X})=\delta_{2,0,0}^{(0)}=X_2$.
\\~\\\textbf{Case II:} $\rho \sigma_1<\sigma_2$ \vspace*{2mm}

In this case $\alpha_{2,0}=\beta_0>0,\; \alpha_{2,\infty}=\infty,\;\psi_2(z)=\frac{\sigma_2(\sigma_2-\rho\sigma_1)}{\tau^2}z$ is increasing in $z\in (-\infty,0)$ and, for every fixed $\lambda\geq 0$, $k_2(z,\lambda)=\frac{\sigma_2(\sigma_2-\rho\sigma_1)}{\tau^2}\frac{z-\lambda}{z}$ is increasing in $z\in (-\infty,0)$. Using Theorem 3.2.1 (a), we conclude that the estimators $\{\delta_{2,0,\alpha}^{(0)}:\alpha\in [\beta_0,\infty)\}$ are admissible within the class $\mathcal{D}_{2,0}=\{\delta_{2,0,\alpha}^{(0)}:-\infty<\alpha<\infty\}$ of isotonic regression estimators of $\theta_2$. Moreover, the estimators $\{\delta_{2,0,\alpha}^{(0)}:\alpha<\beta_0\}$ are inadmissible for estimating $\theta_2$ and, for $\infty< \alpha_1<\alpha_2\leq \beta_0$, the estimator $\delta_{2,0,\alpha_2}^{(0)}$ dominates the estimator $\delta_{2,0,\alpha_1}^{(0)}$. In particular the RMLE $\delta_{2,0,\beta_0}^{(0)}$ is admissible  within the class $\mathcal{D}_{2,0}$ of isotonic regression estimators of $\theta_2$ and it dominates the BLEE $\delta_{2,0}(\underline{X})=\delta_{2,0,0}^{(0)}=X_2$. For $\rho=0$, this result is also reported in Kubokawa and Saleh (\citeyear{MR1370413}).	
\\~\\\textbf{Case III:} $\rho \sigma_1=\sigma_2$ \vspace*{2mm}

In this case $\alpha_{2,0}=\beta_0=\alpha_{2,\infty}=1$ and the BLEE $\delta_{2,0}(\underline{X})=\delta_{2,0,0}^{(0)}(\underline{X})=X_2$ is also the RMLE. Also, $\psi_2(z)=0=k_2(z,\lambda),\; z\in(-\infty,0),\;\lambda\geq 0$. Using Theorem 3.2.1, it follows that the BLEE/RMLE $\delta_{2,0}(\underline{X})=\delta_{2,0,0}^{(0)}(\underline{X})=X_2$ is the only admissible estimator within the class $\mathcal{D}_{2,0}$ of isotonic regression estimators of $\theta_2$. Any other isotonic regression estimator in $\mathcal{D}_{2,0}$ is dominated by the BLEE/RMLE $\delta_{2,0}(\underline{X})=\delta_{2,0,0}^{(0)}(\underline{X})=X_2$.
\vspace*{2mm}

\noindent
\textbf{Example 3.2.3.} For the exponential probability model considered in Example 3.1.3, the BLEE of $\theta_2$ is $\delta_{2,0}(\underline{X})=\delta_{2,\sigma_2,0}^{(\beta)}(\underline{X})=X_2-\sigma_2$. Moreover
$$f_Z(z)= \frac{e^{\frac{z}{\sigma_1}}\,e^{-\left(\frac{1}{\sigma_1}+\frac{1}{\sigma_2}\right) \max\{0,z\}}}{\sigma_1+\sigma_2} =\begin{cases}\frac{e^{\frac{z}{\sigma_1}}}{\sigma_1+\sigma_2},&\text{ if } z<0\\~\\
	\frac{e^{\frac{-z}{\sigma_2}}}{\sigma_1+\sigma_2}, &\text{ if } z\geq 0 \end{cases} $$		
is log-concave on $\Re$,
$\psi_2(z)=E[Z_2|Z=z]=z+\psi_1(z)=\max\{0,z\}+\frac{\sigma_1 \sigma_2}{\sigma_1+\sigma_2}$ is increasing on $\Re$ and, for any $\lambda \geq 0$,
$$k_2(z,\lambda)=\frac{\psi_2(z+\sigma_2-\beta-\lambda)-\sigma_2}{z}=\begin{cases}\frac{-\sigma_2^2}{(\sigma_1+\sigma_2)z},& \text{ if} \;\; z<\lambda+\beta-\sigma_2\\
	1-\left(\lambda+\beta-\frac{\sigma_1 \sigma_2}{\sigma_1+\sigma_2}\right)\frac{1}{z}, &\text{ if} \;\; z\geq \lambda+\beta-\sigma_2
\end{cases}$$
is increasing in $z\in(-\infty,0)$, provided $\beta\geq \frac{\sigma_1 \sigma_2}{\sigma_1+\sigma_2}$. Since $c_{0,1}=\sigma_1$, the suitable choice of $\beta$ in this case is $\beta=\sigma_1$. For $\beta=\sigma_1$, after some straightforward algebra, one gets

\begin{align*}
	\alpha_{2,0}
	&=\frac{\int_{-\infty}^{0}(\psi_2(z+\sigma_2-\sigma_1)-\sigma_2)\,z\,f_Z(z+\sigma_2-\sigma_1)dz}{\int_{-\infty}^{0}z^2\,f_Z(z+\sigma_2-\sigma_1)dz}\\
	&=\begin{cases}
		\frac{\frac{\sigma_2^2(\sigma_1^2-2\sigma_1\sigma_2-2\sigma_2^2)}{\sigma_1+\sigma_2}e^{\frac{\sigma_1}{\sigma_2}}+\sigma_2^2(\sigma_1+\sigma_2)e}{(\sigma_1+\sigma_2)(\sigma_1^2+\sigma_2^2)e-2\sigma_2^3 e^{\frac{\sigma_1}{\sigma_2}}},&\text{if } \sigma_1\leq \sigma_2\\~\\
		\frac{\sigma_2^2}{2\sigma_1(\sigma_1+\sigma_2)},&\text{if } \sigma_1> \sigma_2 \end{cases}\end{align*} 
\begin{align*}	\text{and} \qquad \alpha_{2,\infty}&=\lim_{\lambda\to \infty} \frac{\int_{-\infty}^{0} (\psi_2(z+\sigma_2-\sigma_1-\lambda)-\sigma_2)\,z\, f_Z(z+\sigma_2-\sigma_1-\lambda)dz}{\int_{-\infty}^{0} z^2 f_Z(z+\sigma_2-\sigma_1-\lambda)dz}\\
	&=\frac{\sigma_2^2}{2\sigma_1(\sigma_1+\sigma_2)}.
\end{align*} 
Note that, for $\sigma_1>\sigma_2$,
$$\alpha_2(\lambda)=\alpha_{2,0}=\alpha_{2,\infty}=\frac{\sigma_2^2}{2\sigma_1(\sigma_1+\sigma_2)},\;\forall\,\lambda\geq 0.$$
It is easy to verify that $\alpha_{2,0}>0$, for any $\sigma_1>0$ and $\sigma_2>0$.
Using Theorem 3.2.1 (a), it follows that the estimator estimators $\{\delta_{2,\sigma_2,\alpha}^{(\sigma_1)}: \alpha\in[\alpha_{2,0},\alpha_{2,\infty}] \}$ are admissible within the class $\mathcal{D}_{2,\sigma_1}=\{\delta_{2,\sigma_2,\alpha}^{(\sigma_1)}: -\infty<\alpha<\infty \}$ of isotonic regression estimators of $\theta_2$; here $\beta=\sigma_1$. Moreover, for any $-\infty<\alpha_1<\alpha_2\leq \alpha_{2,0}$ or $\alpha_{2,\infty}\leq \alpha_2<\alpha_1<\infty$, the estimator $\delta_{2,\sigma_2,\alpha_2}^{(\sigma_1)}$ dominates the estimator $\delta_{2,\sigma_2,\alpha_1}^{(\sigma_1)}$. In particular the BLEE $\delta_{2,\sigma_2}(\underline{X})=\delta_{2,\sigma_2,0}^{(\sigma_1)}(\underline{X})=X_2-\sigma_2$ is inadmissible for estimating $\theta_2$ and is dominated by the isotonic regression estimator
\begin{align*}
	\delta_{2,\sigma_2,\alpha_{2,0}}^{(\sigma_1)}(\underline{X})&=\begin{cases}
		X_2-\sigma_2,&\text{ if} \;\; X_1-\sigma_1\leq X_2-\sigma_2\\
		\alpha_{2,0}(X_1-\sigma_1)+(1-\alpha_{2,0})(X_2-\sigma_2), &\text{ if} \; \; X_1-\sigma_1> X_2-\sigma_2
	\end{cases}\\
	& =\max\{X_2-\sigma_2, \alpha_{2,0}(X_1-\sigma_1)+(1-\alpha_{2,0})(X_2-\sigma_2)\}.
\end{align*}
Interestingly, for $\sigma_1>\sigma_2$, the only admissible estimator within the class $\mathcal{D}_{2,\sigma_1}$ of isotonic regression estimators is $\delta_{2,\sigma_2,\alpha_0}^{(\sigma_1)}$, where $\alpha_0=\frac{\sigma_2^2}{2\sigma_1(\sigma_1+\sigma_2)}$, and it dominates the BLEE 
$\delta_{2,\sigma_2}(\underline{X})=X_2-\sigma_2$. Similar results have also been obtained by Misra and Singh (\citeyear{MR1366828}).

\subsection{\textbf{Simulation Study For Estimation of Location Parameter $\theta_1$}}
\label{sec:3.3}
\noindent
\vspace*{2mm}

Under the squared error loss function, in Example 3.1.3, we considered estimation of the smaller location parameter $\theta_1$ of two independent exponential distributions having unknown order restricted location parameters (i.e., $\theta_1\leq \theta_2$) and known scale parameters ($\sigma_1>0$ and $\sigma_2>0$). We have shown that the isotonic regression estimator (IRE) $\min\{X_1-\sigma_1,X_2-\frac{\sigma_1\sigma_2}{\sigma_1+\sigma_2}\}$, based on $(X_1-\sigma_1,X_2-\frac{\sigma_1\sigma_2}{\sigma_1+\sigma_2})$ with $\alpha=0$, dominates the BLEE $\delta_{0,1}(\underline{X})=X_1-\sigma_1$. To further evaluate the performances of various estimators under the squared error loss function, in this section, we compare the risk performances of the BLEE $\delta_{0,1}(\underline{X})=X_1-\sigma_1$, the Restricted MLE $\min\{X_1,X_2\}$, the IRE $\min\{X_1-\sigma_1,X_2-\sigma_2\}$, based on the BLEEs $(X_1-\sigma_1,X_2-\sigma_2)$ with $\alpha=0$, and the IRE $\min\{X_1-\sigma_1,X_2-\frac{\sigma_1\sigma_2}{\sigma_1+\sigma_2}\}$, based on $(X_1-\sigma_1,X_2-\frac{\sigma_1\sigma_2}{\sigma_1+\sigma_2})$ with $\alpha=0$, numerically, through Monte Carlo simulations. For simulations, we generated 10000 samples of size 1 each from relevant exponential distributions and computed the simulated risks of the BLEE, the Restricted MLE, the Isotonic regression Estimator (IRE) based on $(X_1-\sigma_1,X_2-\sigma_2)$ with $\alpha=0$, and the IRE based on $(X_1-\sigma_1,X_2-\frac{\sigma_1\sigma_2}{\sigma_1+\sigma_2})$ with $\alpha=0$. 
\vspace*{2mm}

\begin{figure}[h!]
	\begin{subfigure}{0.48\textwidth}
		\centering
		\includegraphics[width=72mm,scale=1.2]{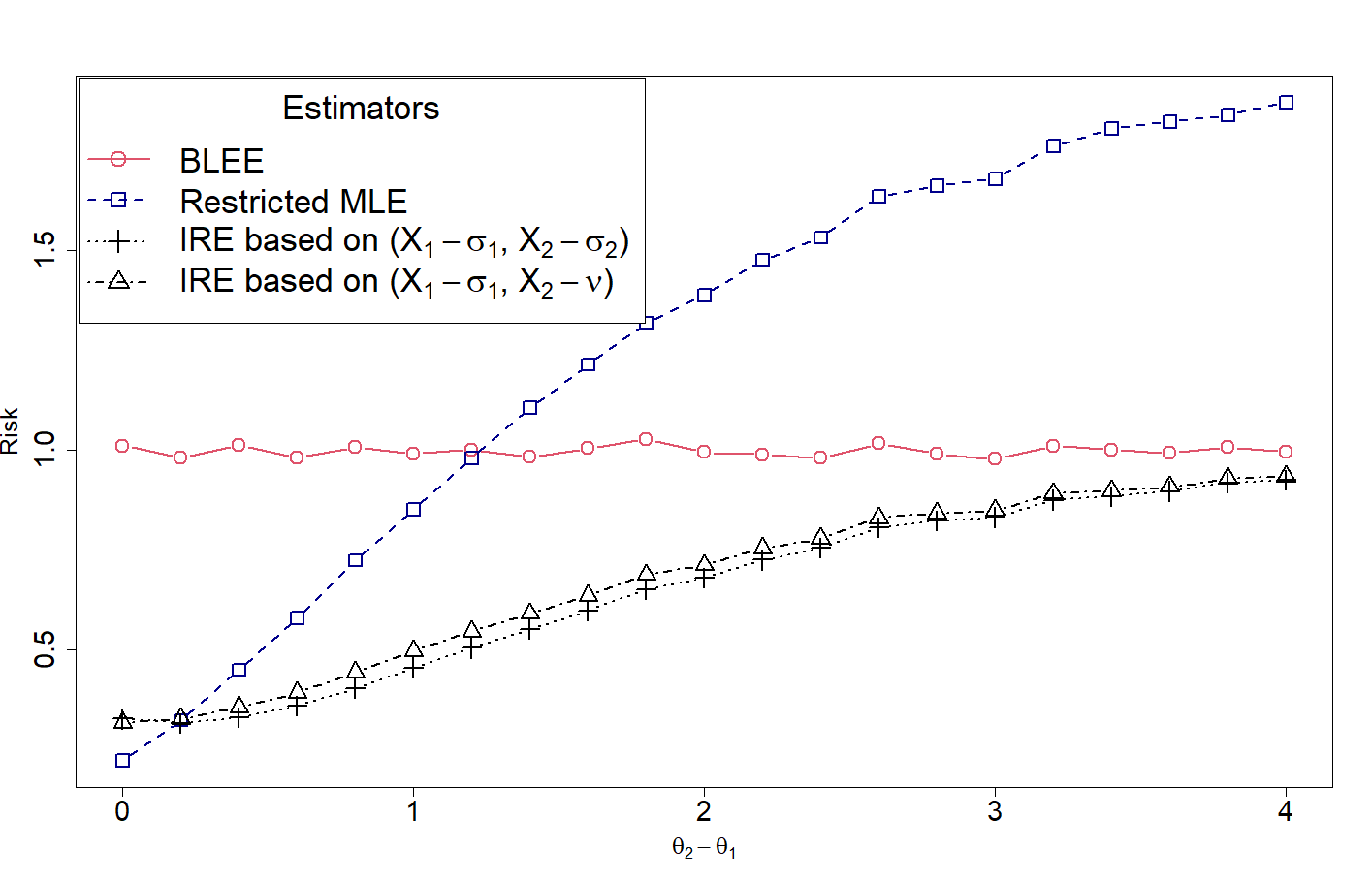} 
		\caption{$\sigma_1=1$ and $\sigma_2=0.5$.} 
		\label{fig1:a} 
	\end{subfigure}
	\begin{subfigure}{0.48\textwidth}
		\centering
		\includegraphics[width=72mm,scale=1.2]{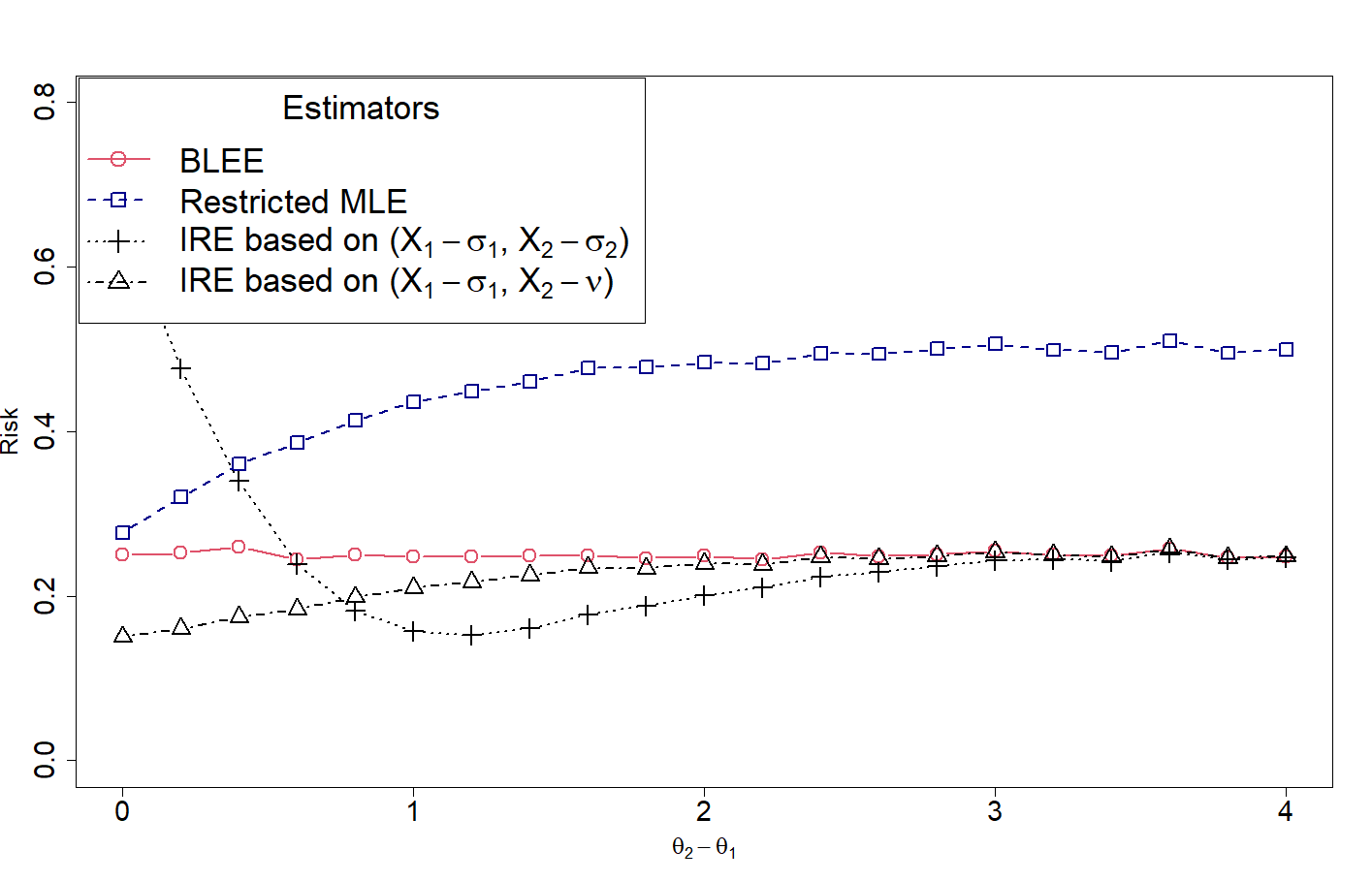} 
		\caption{$\sigma_1=0.5$ and $\sigma_2=1.5$.} 
		\label{fig1:b} 
	\end{subfigure}
	\\	\begin{subfigure}{0.48\textwidth}
		\centering
		
		\includegraphics[width=72mm,scale=1.2]{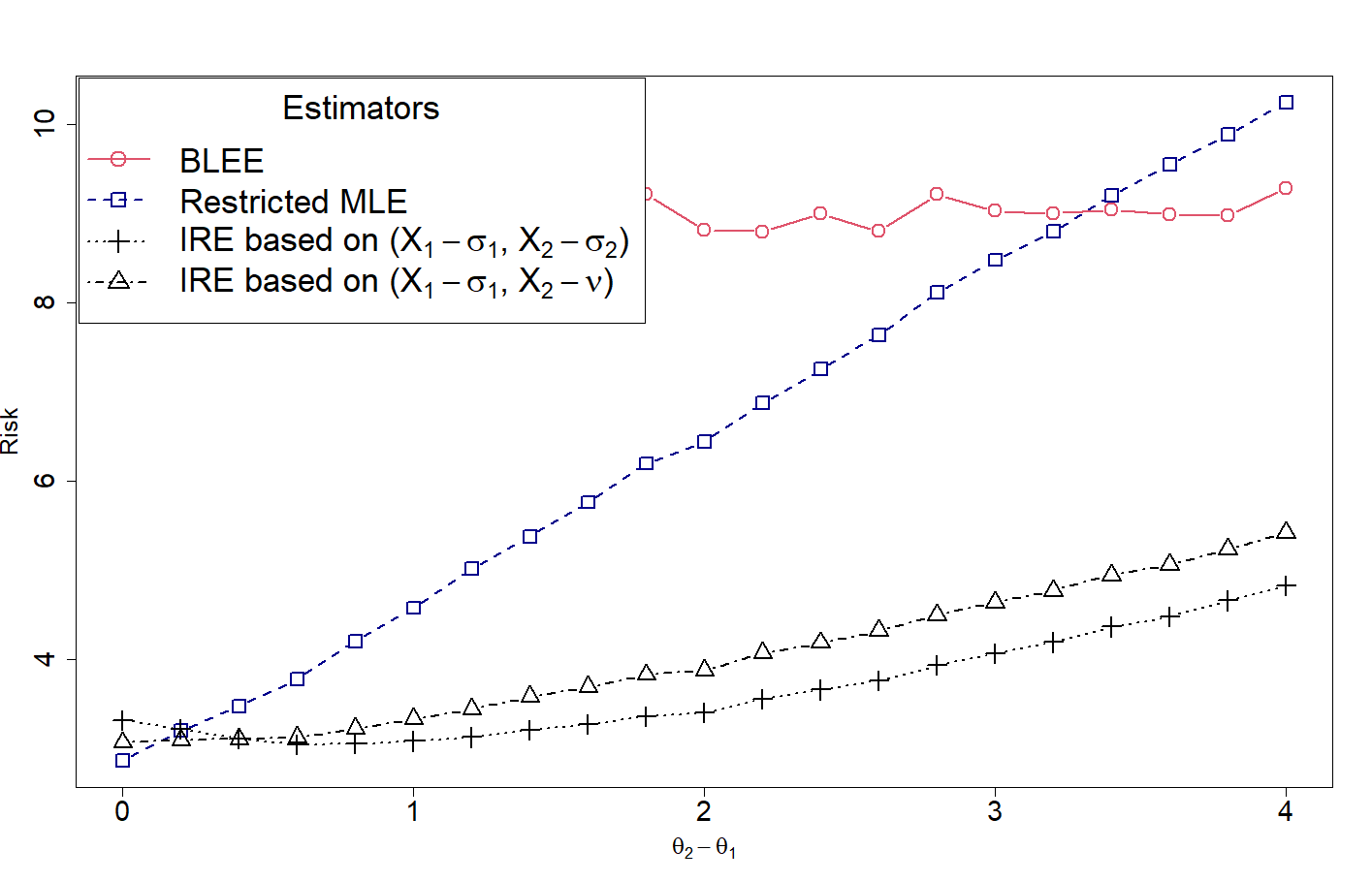} 
		
		\caption{$\sigma_1=3$ and $\sigma_2=2$.} 
		\label{fig1:c} 
	\end{subfigure}
	\begin{subfigure}{0.48\textwidth}
		\centering
		\includegraphics[width=72mm,scale=1.2]{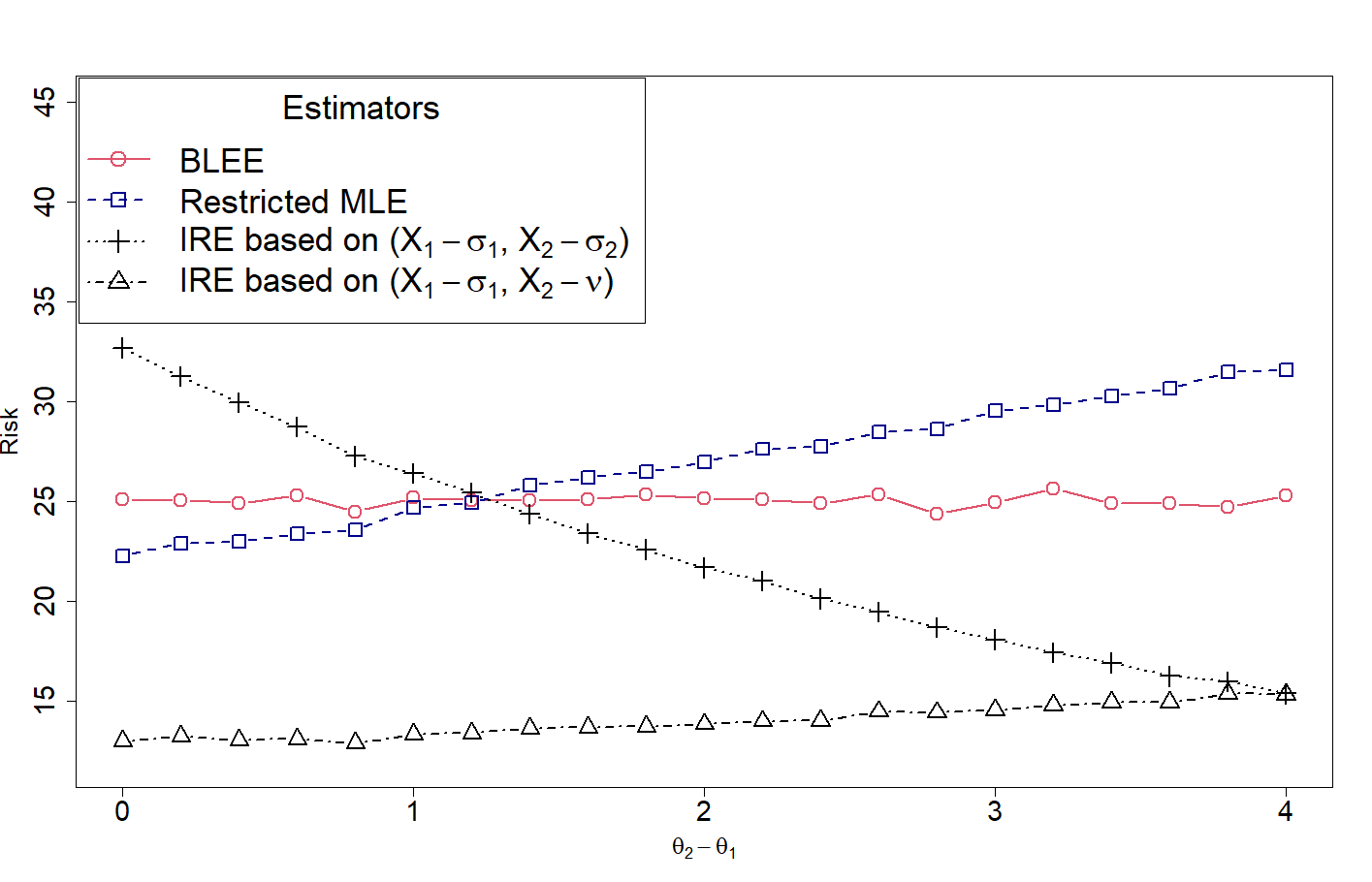} 
		\caption{ $\sigma_1=5$ and $\sigma_2=10$.} 
		\label{fig1:d}  
	\end{subfigure}
	\\	\begin{subfigure}{0.48\textwidth}
		\centering
		
		\includegraphics[width=72mm,scale=1.2]{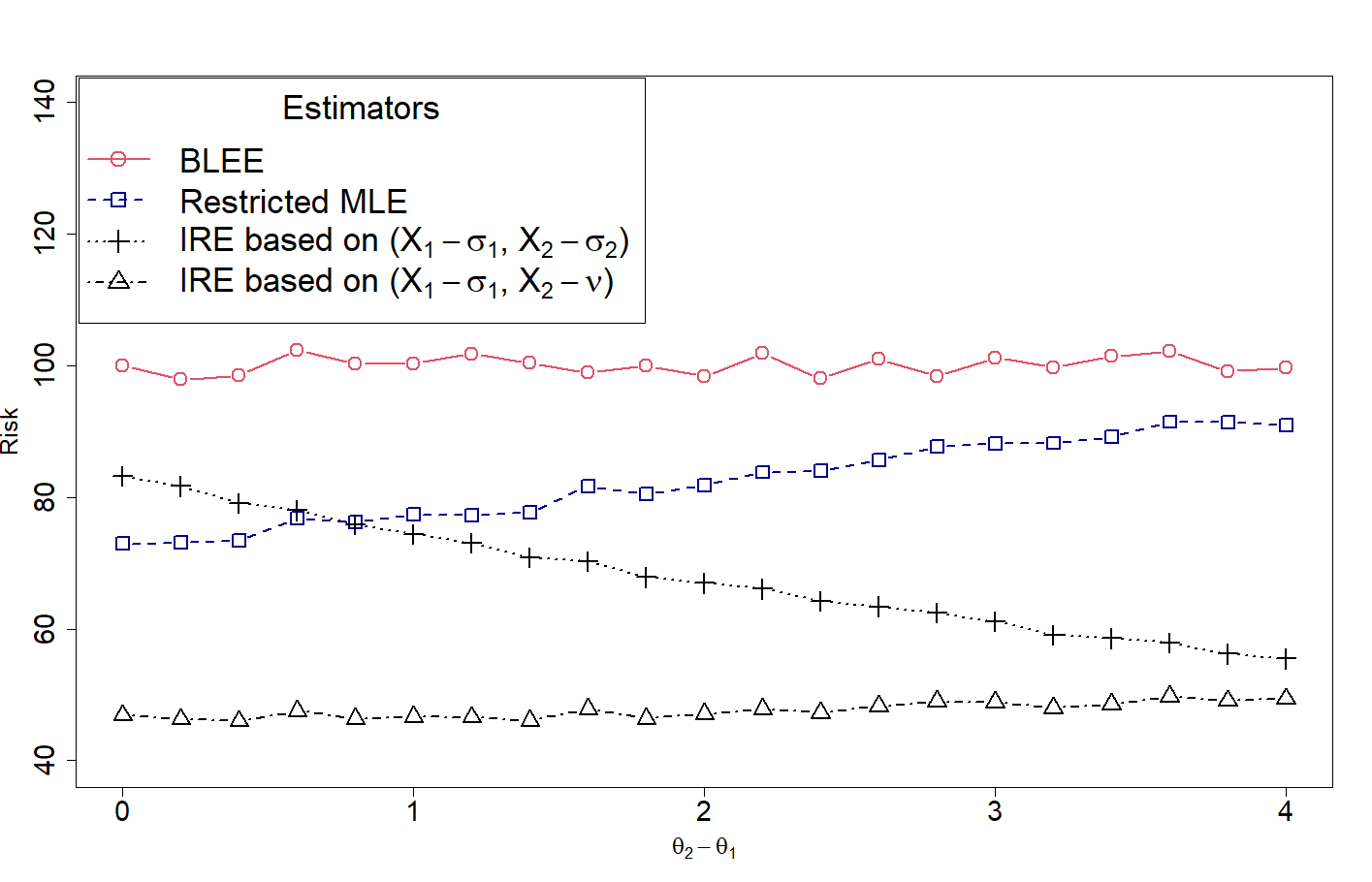} 
		
		\caption{$\sigma_1=10$ and $\sigma_2=15$.} 
		\label{fig1:e} 
	\end{subfigure}
	\begin{subfigure}{0.48\textwidth}
		\centering
		\includegraphics[width=72mm,scale=1.2]{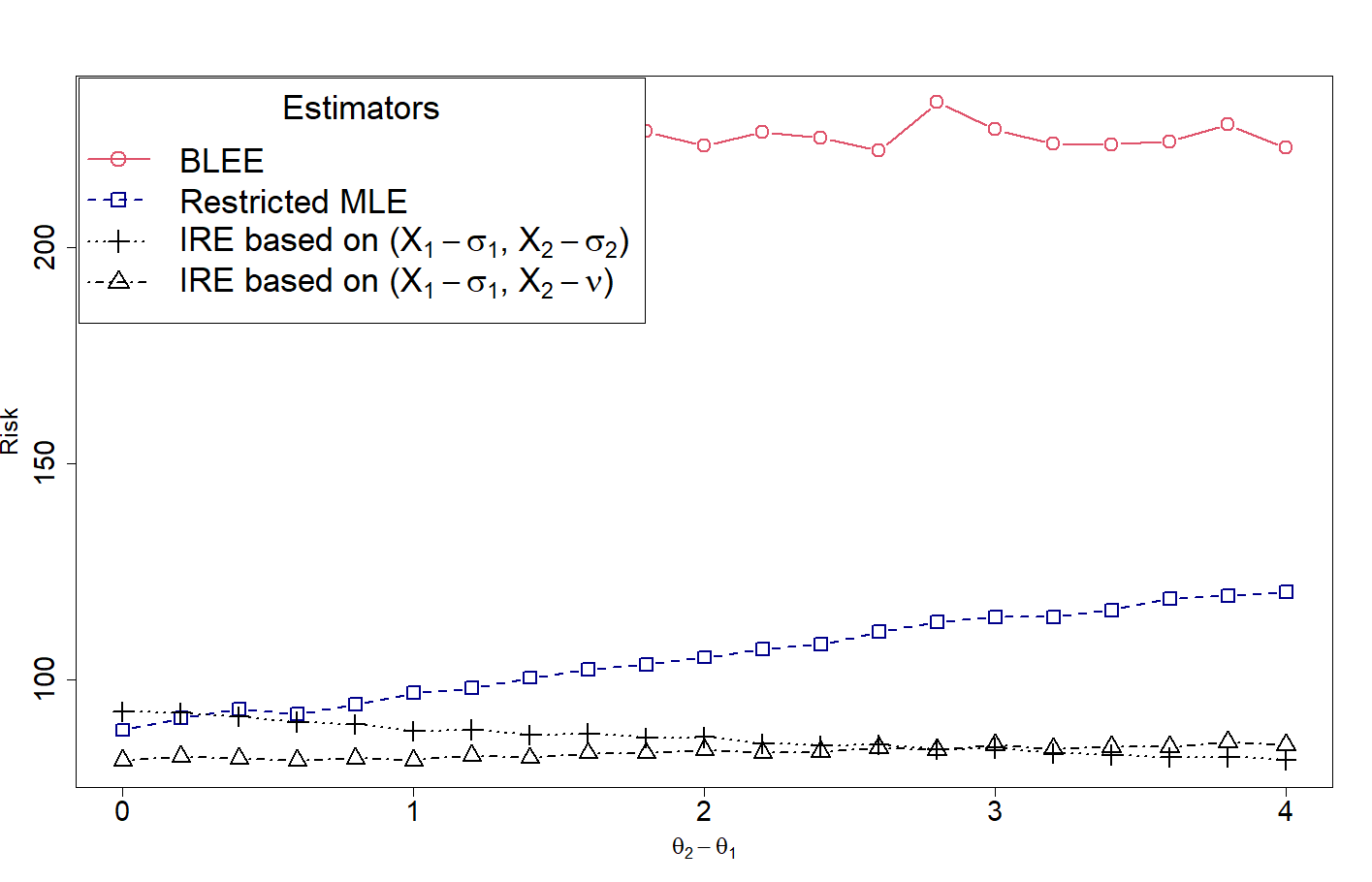} 
		\caption{ $\sigma_1=15$ and $\sigma_2=12$.} 
		\label{fig1:f}  
	\end{subfigure}
	\caption{Risk plots of the BLEE $X_1-\sigma_1$, the Restricted MLE $\min\{X_1,X_2\}$, the IRE based on $(X_1-\sigma_1,X_2-\sigma_2)$ with $\alpha=0$, and the IRE based on $(X_1-\sigma_1,X_2-\nu)$ with $\alpha=0$, where $\nu=\frac{\sigma_1\sigma_2}{\sigma_1+\sigma_2}$, against the values of, $\lambda=\theta_2-\theta_1$.}
	\label{fig1}
\end{figure}

The simulated values of risks of various estimators are plotted in Figure \ref{fig1}. The following observations are evident from Figure \ref{fig1}:\vspace*{1.5mm}

\noindent
(i) The risk function values of the IRE based on $(X_1-\sigma_1,X_2-\frac{\sigma_1\sigma_2}{\sigma_1+\sigma_2})$ with $\alpha=0$ estimator is nowhere larger than the risk function values of the BLEE, which is in 
conformity with theoretical findings of Example 3.1.3.\vspace*{1.5mm}

\noindent
(ii) Interestingly, here IRE based on the BLEEs $(X_1-\sigma_1,X_2-\sigma_2)$ with $\alpha=0$ does not always dominates the BLEE $X_1-\sigma_1$.\vspace*{1.5mm}

\noindent
\\(iii) There is no clear cut winner between various estimators. But the IRE based on $(X_1-\sigma_1,X_2-\frac{\sigma_1\sigma_2}{\sigma_1+\sigma_2})$ with $\alpha=0$, performs reasonably well.
\FloatBarrier
\FloatBarrier

\section{\textbf{Isotonic Regression Estimators For Component-wise Estimation of Order Restricted Scale Parameters}}
\label{sec:4}
Let $\underline{X}=(X_1,X_2)$ be a random vector having the Lebesgue pdf belonging to the scale family \eqref{eq:1.2}. We assume that the distributional support of $\underline{X}=(X_1,X_2)$ is a subset of $\Re_{++}^{2}$.
Consider estimation of $\theta_i$ under the scaled squared error loss (SSEL) function \eqref{eq:1.4}, when it is known apriori that $\underline{\theta}\in \Theta_0=\{(x_1,x_2)\in \Re_{++}:x_1\leq x_2\}$. 
Let $c_{0,i}$ be as defined in (1.6), so that $\delta_{i,c_{0,i}}(\underline{X})=c_{0,i} X_i$ is the BSEE of $\theta_i,\;i=1,2$. For suitably chosen $\nu >0$ and $\beta>0$, let $\mathcal{D}_{1,\nu}$ and $\mathcal{D}_{2,\beta}$ be the classes of isotonic regression estimators of $\theta_1$ and $\theta_2$ defined by (1.9) and (1.10), respectively, where $\delta_{2,\nu}(\underline{X})=\nu X_2$ and $\delta_{1,\beta}(\underline{X})=\beta X_1$. Here $\nu$ and $\beta$ will be chosen in a manner such that SEEs $\delta_{2,\nu}(\underline{X})=\nu  X_2$ and $\delta_{1,\beta}(\underline{X})=\beta  X_1$ are close to the BSEEs $\delta_{2,c_{0,2}}(\underline{X})=c_{0,2} X_2$ and $\delta_{1,c_{0,1}}(\underline{X})=c_{0,1} X_1$, respectively, and dominance of $\delta_{1,c_{0,1,\alpha}}^{(\nu)}(\underline{X})$ and $\delta_{2,c_{0,2},\alpha}^{(\beta)}(\underline{X})$ over the BSEEs $\delta_{1,c_{0,1}}$ and $\delta_{2,c_{0,2}}$, respectively, is ensured for some $\alpha$'s. Thus, to begin with, we take $\nu$ and $\beta$ to be fixed positive numbers and consider classes  $\mathcal{D}_{1,\nu}$ and $\mathcal{D}_{2,\beta}$ of isotonic regression estimators of $\theta_1$ and $\theta_2$, respectively. We aim to characterize admissible estimators within classes $\mathcal{D}_{1,\nu}$ and $\mathcal{D}_{2,\beta}$ of isotonic regression estimators of $\theta_1$ and $\theta_2$, respectively, and to find estimators in these classes that dominate the BSEEs $\delta_{1,c_{0,1}}(\underline{X})=c_{0,1} X_1$ and $\delta_{1,c_{0,2}}(\underline{X})=c_{0,2} X_2$. The following subsection deals with estimation of the smaller scale parameter $\theta_1$.

\subsection{\textbf{Isotonic Regression Estimators of the smaller Scale Parameter $\theta_1$}}
\label{sec:4.1}
\noindent
\vspace*{2mm}

Let $\nu$ be a fixed positive constant, to be suitably chosen as described above. Consider the class $\mathcal{D}_{1,\nu}$ of isotonic regression estimators defined by (1.7) and (1.9), where $\delta_{1,c_{0,1}}(\underline{X})=c_{0,1}X_1$ and $\delta_{2,\nu}(\underline{X})=\nu X_2$. Define $D=\frac{X_2}{X_1}$, $\lambda=\frac{\theta_2}{\theta_1}$ ($\lambda \geq 1$), $Z_1=\frac{X_1}{\theta_1}$, $Z_2=\frac{X_2}{\theta_2}$, $Z=\frac{Z_2}{Z_1}$ and $\underline{Z}=(Z_1,Z_2)$. For estimating $\theta_1$ under the SEL function $L_2$, defined by (1.4), the risk function of the estimator $\delta_{1,c_{0,1},\alpha}^{(\nu)}\in \mathcal{D}_{1,\nu}$,\; $-\infty<\alpha<\infty$, is given by
\begin{align*}
	R_1(\underline{\theta},\delta_{1,c_{0,1},\alpha}^{(\nu)})&=E_{\underline{\theta}}\left[\left(\frac{\delta_{1,c_{0,1},\alpha}^{(\nu)}(\underline{X})}{\theta_1}-1\right)^2\right]\\
	&=E_{\underline{\theta}}[(c_{0,1}Z_1-1)^2 I_{(\frac{c_{0,1}}{\lambda \nu},\infty)}(Z)] \\
	&\qquad\qquad	+ E_{\underline{\theta}}[(\alpha c_{0,1} Z_1+(1-\alpha)\nu\lambda Z_2-1)^2 I_{(0,\frac{c_{0,1}}{\lambda \nu})}(Z)].
\end{align*}
For any fixed $\underline{\theta}\in\Theta_0$ (or $\lambda \geq 1)$, $R_1(\underline{\theta},\delta_{1,c_{0,1},\alpha}^{(\nu)})$ is minimized at $\alpha= \alpha_1(\lambda)$, where
\begin{align}
	\alpha_1(\lambda)=&\frac{E_{\underline{\theta}}[(1-\nu \lambda Z_2)(c_{0,1}Z_1-\nu \lambda Z_2) I_{(0,\frac{c_{0,1}}{\lambda \nu})}(Z)]}{E_{\underline{\theta}}[(c_{0,1}Z_1-\nu \lambda Z_2)^2 I_{(0,\frac{c_{0,1}}{\lambda \nu})}(Z)]}
	=1+\alpha_1^{*}(\lambda),\;\;\lambda \geq 1, \label{eq:4.1}\\ \nonumber \\ \nonumber
	\text{and}\quad \alpha_1^{*}(\lambda)=&\frac{E_{\underline{\theta}}[(c_{0,1}Z_1-\nu \lambda Z_2)(1-c_{0,1}Z_1) I_{(0,\frac{c_{0,1}}{\lambda \nu})}(Z)]}{E_{\underline{\theta}}[(c_{0,1}Z_1-\nu \lambda Z_2)^2 I_{(0,\frac{c_{0,1}}{\lambda \nu})}(Z)]}\nonumber\\
	=&\frac{1}{c_{0,1}}\frac{E_{\underline{\theta}}[(1-\nu \frac{\lambda Z}{c_{0,1}})(\psi_1(Z)-c_{0,1}\psi_2(Z)) I_{(0,\frac{c_{0,1}}{\lambda \nu})}(Z)]}{E_{\underline{\theta}}[(1-\frac{\nu \lambda Z}{c_{0,1}})^2 I_{(0,\frac{c_{0,1}}{\lambda \nu})}(Z)]}\nonumber\\
	=& \frac{1}{c_{0,1}} \frac{\int_{0}^{1} (1-z) (\psi\left(\frac{c_{0,1} z}{\nu \lambda}\right)-c_{0,1})\, \psi_2\left(\frac{c_{0,1} z}{\nu \lambda}\right)\,f_Z\left(\frac{c_{0,1} z}{\nu \lambda}\right) dz}{\int_{0}^{1}  (1-z)^2 \,\psi_2\left(\frac{c_{0,1} z}{\nu \lambda}\right)\,f_Z\left(\frac{c_{0,1} z}{\nu \lambda}\right) dz},\;\; \lambda\geq 1;\label{eq:4.2}
\end{align}
here $f_Z(z)=\int_{0}^{\infty} s f(s,sz)ds,\;0<z<\infty$ is the pdf of $Z$, $\psi_r(z)=E\left[Z^{^r}_1|Z=z\right]=\frac{\int_{0}^{\infty}s^{r+1} f(s,sz)ds}{\int_{0}^{\infty}s f(s,sz)ds},\; z\in S_Z,\;r=1,2$, and $\psi(z)=\frac{\psi_1(z)}{\psi_2(z)},\;z\in S_Z$.
\\~\\Then,
$
\alpha_1^{*}(\lambda)
=\frac{1}{c_{0,1}}E_{\lambda}[k_1(S_{1,\lambda},\lambda)],\, \lambda \geq 1,
$
where 
$
k_1(z,\lambda)=\frac{\psi\left(\frac{c_{0,1} z}{\nu \lambda}\right)-c_{0,1}}{1-z},\; 0<z<1,\; \lambda \geq 1,
$
and $S_{1,\lambda}$ is a r.v. having the pdf
\newpage
\begin{equation*}
	h_{1,\lambda}(t)=\begin{cases}
		\frac{(1-t)^2 \,\psi_2\left(\frac{c_{0,1} t}{\nu \lambda}\right)\,f_Z\left(\frac{c_{0,1} t}{\nu \lambda}\right) }{\int_{0}^{1}(1-s)^2 \,\psi_2\left(\frac{c_{0,1} s}{\nu \lambda}\right)\,f_Z\left(\frac{c_{0,1} s}{\nu \lambda}\right) ds},& \text{if}\;\; 0<t<1\\
		0, &\text{otherwise}
	\end{cases},\;\;\lambda \geq 1.
\end{equation*}
Clearly if, for every fixed $\theta\in(0,1)$, $\frac{\psi_2(\theta t)}{\psi_2(t)}$ and $\frac{f_Z(\theta t)}{f_Z(t)}$ are increasing (decreasing) functions of t on $(0,\frac{c_{0,1}}{\nu})$, then $S_{1,\lambda_1}\leq_{lr}S_{1,\lambda_2}$ ($S_{1,\lambda_2}\leq_{lr}S_{1,\lambda_1}$) and, consequently, $S_{1,\lambda_1}\leq_{st}S_{1,\lambda_2}$ ($S_{1,\lambda_2}\leq_{st}S_{1,\lambda_1}$), whenever $1\leq \lambda_1<\lambda_2<\infty$. \vspace*{2mm}

\noindent
The following lemma, whose proof is immediate from Proposition 2.1, will be useful in arriving at the main result of this subsection.
\\~\\ \textbf{Lemma 4.1.1 (a)} Suppose that, for every fixed $\theta\in(0,1)$, $\frac{\psi_2(\theta t)}{\psi_2(t)}$ and $\frac{f_Z(\theta t)}{f_Z(t)}$ are increasing (decreasing) functions of t on $(0,\frac{c_{0,1}}{\nu})$. Also, suppose that $\psi(z)$ is decreasing on $(0,\frac{c_{0,1}}{\nu})$ and, for every $\lambda \geq 1$, $k_1(z,\lambda)=\frac{\psi\left(\frac{c_{0,1} z}{\nu \lambda}\right)-c_{0,1}}{1-z}$ is increasing (decreasing) in $z\in(0,1)$. Then $\alpha_1^{*}(\lambda)$ (and hence $\alpha_1(\lambda)$), defined by \eqref{eq:4.2} (\eqref{eq:4.1}), is an increasing function of $\lambda\in[1,\infty)$,
$\inf_{\lambda \geq 1} \alpha_1(\lambda)= 1+\alpha_1^{*}(1)=\alpha_{1,1}, \text{ say, and }
\sup_{\lambda \geq 1} \alpha_1(\lambda)= 1+\lim_{\lambda \to \infty}\alpha_1^{*}(\lambda)=\alpha_{1,\infty}  \text{, say}.$
\\~\\\textbf{(b)}  Suppose that, for every fixed $\theta\in(0,1)$, $\frac{\psi_2(\theta t)}{\psi_2(t)}$ and $\frac{f_Z(\theta t)}{f_Z(t)}$ are increasing (decreasing) functions of t on $(0,\frac{c_{0,1}}{\nu})$. Also, suppose that $\psi(z)$ is increasing on $(0,\frac{c_{0,1}}{\nu})$ and, for every $\lambda \geq 1$, $k_1(z,\lambda)=\frac{\psi\left(\frac{c_{0,1} z}{\nu \lambda}\right)-c_{0,1}}{1-z}$ is decreasing (increasing) in $z\in(0,1)$. Then $\alpha_1^{*}(\lambda)$ (and hence $\alpha_1(\lambda)$) is a decreasing function of $\lambda\in[1,\infty)$,
$\inf_{\lambda \geq 1} \alpha_1(\lambda) = 1+\lim_{\lambda \to \infty}\alpha_1^{*}(\lambda)=\alpha_{1,\infty}, \text{ say, and }
\sup_{\lambda \geq 1} \alpha_1(\lambda)= 1+\alpha_1^{*}(1)=\alpha_{1,1} \text{, say}.$
\vspace*{2mm}

\noindent
Now we state the main result of this subsection. The proof of the theorem, being similar to the proof of Theorem 3.1.1, is omitted.
\\~\\ \textbf{Theorem 4.1.1 (a)} Suppose that assumptions of Lemma 4.1.1 (a) hold. Then the estimators that are admissible in the class $\mathcal{D}_{1,\nu}$ are $\{\delta_{1,c_{0,1},\alpha}^{(\nu)}:\alpha\in[\alpha_{1,1},\alpha_{1,\infty}]\}$. Moreover, for $-\infty<\alpha_1<\alpha_2\leq \alpha_{1,1}$ or $\alpha_{1,\infty}\leq \alpha_2<\alpha_1$, the estimator $\delta_{1,c_{0,1},\alpha_2}^{(\nu)}$ dominates the estimator $\delta_{1,c_{0,1},\alpha_1}^{(\nu)}$, for any $\underline{\theta}\in\Theta_0$.
\\~\\ \textbf{(b)} Suppose that assumptions of Lemma 4.1.1 (b) hold. Then the estimators that are admissible in the class $\mathcal{D}_{1,\nu}$ are $\{\delta_{1,c_{0,1},\alpha}^{(\nu)}:\alpha\in[\alpha_{1,\infty},\alpha_{1,1}]\}$. Moreover, for $-\infty<\alpha_1<\alpha_2\leq \alpha_{1,\infty}$ or $\alpha_{1,1}\leq \alpha_2<\alpha_1<\infty$, the estimator $\delta_{1,c_{0,1},\alpha_2}^{(\nu)}$ dominates the estimator $\delta_{1,c_{0,1},\alpha_1}^{(\nu)}$, for any $\underline{\theta}\in\Theta_0$.
\vspace*{2mm}

\noindent Now we illustrate some applications of Theorem 4.1.1.
\\~\\ \textbf{Example 4.1.1.} Let $f(z_1,z_2)=f_1(z_1) f_2(z_2),\;\; \underline{z}=(z_1,z_2)\in \Re^2,$ where, for positive constants $\alpha_1$ and $\alpha_2$,
$$f_i(z)=\begin{cases}  \frac{e^{-z} z^{\alpha_i-1}}{\Gamma{\alpha_i}}, &\text{ if} \;\; z>0 \\
	0, & \text{ otherwise}  \end{cases},\;\; i=1,2.  $$
Here $c_{0,i}=\frac{E[Z_i]}{E[Z_i^2]}=\frac{1}{\alpha_i+1},\;i=1,2,$ and the BSEE of $\theta_1$ is $\delta_{1,0}(\underline{X})=\delta_{1,0,1}^{(\nu)}(\underline{X})=\frac{X_1}{\alpha_1+1}$. The pdf of $Z=\frac{Z_2}{Z_1}$ is 
$$f_Z(z)= \begin{cases} \frac{\Gamma{(\alpha_1+\alpha_2)}}{\Gamma{\alpha_1} \Gamma{\alpha_2}}\frac{z^{\alpha_2-1}}{(1+z)^{\alpha_1+\alpha_2}}, &\text{ if} \;\; 0<z<\infty \\
	0, & \text{ otherwise}  \end{cases}.  $$
One can easily check that, for every fixed $\theta\in(0,1)$, $\frac{f_Z(\theta t)}{f_Z(t)}$ is increasing in $t\in(0,\infty)$, the conditional pdf of $Z_1$ given $Z=z$ ($0<z<\infty$) is
$$ f_{Z_1|Z}(z_1|z)=\begin{cases}
	\frac{ (1+z)^{\alpha_1+\alpha_2}\, z_1^{\alpha_1+\alpha_2-1}\,  e^{-(1+z)z_1}}{\Gamma{(\alpha_1+\alpha_2)}}, &\text{ if} \;\; 0<z_1<\infty \\
	0, & \text{ otherwise}  \end{cases}, $$
$\psi_1(z)=E[Z_1|Z=z]=\frac{\alpha_1+\alpha_2}{1+z},\; z>0,$ and 
$\psi_2(z)=E[Z_1^2|Z=z]=\frac{(\alpha_1+\alpha_2+1)(\alpha_1+\alpha_2)}{(1+z)^2},\; z>0.$
Clearly, for every fixed $\theta\in(0,1)$, $\frac{\psi_2(\theta z)}{\psi_2(z)}=\frac{(1+z)^2}{(1+\theta z)^2}$ is increasing in $z\in (0,\infty)$, and $\psi(z)=\frac{\psi_1(z)}{\psi_2(z)}=\frac{1+z}{\alpha_1+\alpha_2+1}$ is increasing in $z\in (0,\infty)$. For any fixed $\lambda\geq 1$,
$$ k_1(z,\lambda)=\frac{\psi\left(\frac{c_{0,1}z}{\nu \lambda}\right)-c_{0,1}}{1-z}=\frac{z-\nu \lambda \alpha_0}{(\alpha_1+1)(\alpha_1+\alpha_2+1)\nu \lambda(1-z)},\; z>0,$$
is decreasing in $z\in (0,1)$, provided $\nu\geq \frac{1}{\alpha_2}$. Since $c_{0,2}=\frac{1}{\alpha_2+1}$, an appropriate choice of $\nu$ in this case is $\nu=\frac{1}{\alpha_2}$. For $\nu=\frac{1}{\alpha_2}$, it is easy to verify that
\begin{align*}
	\alpha_{1,1}&=1+\alpha_1^{*}(1)\\
	&=1+\frac{(\alpha_1+1)\int_{0}^{1} (1-t) \left(\psi\left(\frac{\alpha_2 t}{\alpha_1+1}\right)-c_{0,1}\right) \psi_2\left(\frac{\alpha_2 t}{\alpha_1+1}\right) f_Z\left(\frac{\alpha_2 t}{\alpha_1+1}\right) dt}{\int_{0}^{1} (1-t)^2 \psi_2\left(\frac{\alpha_2 t}{\alpha_1+1}\right) f_Z\left(\frac{\alpha_2 t}{\alpha_1+1}\right) dt}\\
	&=\frac{\alpha_1+1}{\alpha_1+\alpha_2+1}=\beta_0,\;\;(say),\\~\\
	\text{and}\quad 
	\alpha_{1,\infty}&=1+\lim_{\lambda \to \infty}\alpha_1(\lambda)\\
	&=	1+\lim_{\lambda \to \infty}\frac{(\alpha_1+1)\int_{0}^{1} (1-t) \left(\psi\left(\frac{\alpha_2 t}{\lambda(\alpha_1+1)}\right)-c_{0,1}\right) \psi_2\left(\frac{\alpha_2 t}{\lambda(\alpha_1+1)}\right) f_Z\left(\frac{\alpha_2 t}{\lambda(\alpha_1+1)}\right) dt}{\int_{0}^{1} (1-t)^2 \,\psi_2\left(\frac{\alpha_2 t}{\lambda(\alpha_1+1)}\right) f_Z\left(\frac{\alpha_2 t}{\lambda(\alpha_1+1)}\right) dt}\\
	&=1-\frac{\alpha_2(\alpha_2+2)}{2(\alpha_1+\alpha_2+1)}=\alpha_0,\;\; (say).
\end{align*}
Using Theorem 4.1.1. (b), it follows that the estimators $\{\delta_{1,c_{0,1},\alpha}^{(\nu)}: \alpha\in \left[\alpha_0,\beta_0	\right]\}$ are admissible within the class  $\mathcal{D}_{1,\nu}=\{\delta_{1,c_{0,1},\alpha}^{(\nu)}: -\infty <\alpha < \infty\}$ of isotonic regression estimators of $\theta_1$; here $c_{0,1}=\frac{1}{\alpha_1+1}$ and $\nu=\frac{1}{\alpha_2}$. Moreover, for $-\infty<\alpha_1<\alpha_2\leq \alpha_0$ or $\beta_0\leq \alpha_2<\alpha_1<\infty$, the estimator $\delta_{1,c_{0,1},\alpha_2}^{(\nu)}$ dominates the estimator $\delta_{1,c_{0,1},\alpha_1}^{(\nu)}$. In particular, the BSEE $\delta_{1,c_{0,1}}(\underline{X})=\delta_{1,c_{0,1},1}^{(\nu)}(\underline{X})=\frac{X_1}{\alpha_1+1}$ is inadmissible for estimating $\theta_1$ and is dominated by the isotonic regression estimator
\begin{align*}
	\delta_{1,c_{0,1},\beta_0}^{(\nu)}(\underline{X})=\begin{cases}
		\frac{X_1}{\alpha_1+1}, &\text{ if  } \frac{X_1}{\alpha_1+1}\leq \frac{X_2}{\alpha_2}\\
		\frac{X_1+X_2}{\alpha_1+\alpha_2+1}, &\text{ if  } \frac{X_1}{\alpha_1+1}> \frac{X_2}{\alpha_2}
	\end{cases}
	=\min\bigg\{\frac{X_1}{\alpha_1+1},\frac{X_1+X_2}{\alpha_1+\alpha_2+1}\bigg\}.
\end{align*}
The above estimator is also obtained in Vijayasree et al. (\citeyear{MR1345425}). Note that here using isotonic regression estimator based on $(\frac{X_1}{\alpha_1+1},\frac{X_2}{\alpha_2})$, and not $(\frac{X_1}{\alpha_1+1},\frac{X_2}{\alpha_2+1})$, is successful in providing improvement over the BSEE $\delta_{1,c_{0,1}}(\underline{X})=\frac{X_1}{\alpha_1+1}$.
\\~\\\textbf{Example 4.1.2.} Let $f(z_1,z_2)=f_1(z_1) f_2(z_2),\; \underline{z}=(z_1,z_2)\in \Re_{++}^2,$ where, for positive constants $\alpha_1$ and $\alpha_2$,
$$f_i(z)=\begin{cases}  \alpha_i z^{\alpha_i-1}, &\text{ if} \;\; 0<z<1 \\
	0, & \text{ otherwise}  \end{cases},\;\; i=1,2.  $$
Here $c_{0,i}=\frac{E[Z_i]}{E[Z_i^2]}=\frac{\alpha_i+2}{\alpha_i+1},\;i=1,2,$ and the BSEE of $\theta_1$ is $\delta_{1,c_{0,1}}(\underline{X})=\delta_{1,0,1}^{(\nu)}(\underline{X})=\frac{\alpha_1+2}{\alpha_1+1}X_1$. The pdf of $Z=\frac{Z_2}{Z_1}$ is 
$$f_Z(z)= \frac{\alpha_1 \alpha_2}{\alpha_1+\alpha_2} z^{\alpha_2-1} \left(\min\bigg\{1,\frac{1}{z}\bigg\}\right)^{\alpha_1+\alpha_2},\;\; z>0,  $$
the conditional pdf of $Z_1$ given $Z=z$ ($0<z<\infty$) is
$$ f_{Z_1|Z}(z_1|z)=\begin{cases}
	\frac{(\alpha_1+\alpha_2) \,z_1^{\alpha_1+\alpha_2-1}}{ \left(\min\big\{1,\frac{1}{z}\big\}\right)^{\alpha_1+\alpha_2}}, &\text{ if} \;\; 0<z_1<\min\big\{1,\frac{1}{z}\big\} \\
	0, & \text{ otherwise}  \end{cases} , $$ 
$\psi_1(z)=E[Z_1|Z=z]=\frac{\alpha_1+\alpha_2}{\alpha_1+\alpha_2+1}\min\big\{1,\frac{1}{z}\big\},\, z>0,$ and 
$\psi_2(z)=E[Z_1^2|Z=z]=$ \\ $\frac{\alpha_1+\alpha_2}{\alpha_1+\alpha_2+2}(\min\big\{1,\frac{1}{z}\big\})^2,\, z>0.$
Clearly, for $0<\theta<1$, 
\begin{align*}
	\frac{f_Z(\theta t)}{f_Z(t)}=\begin{cases} \theta^{\alpha_2-1}, &\text{if } \; 0<t< 1\\ \theta^{\alpha_2-1}t, &\text{if }\; 1\leq t< \frac{1}{\theta}\\ \frac{1}{\theta^{\alpha_1+1}}, &\text{if }\; \frac{1}{\theta} \leq t< \infty \end{cases}
	\;\;  \text{ and }\;\;
	\frac{\psi_2(\theta z)}{\psi_2(z)}=\begin{cases} 1, &\text{if } \; 0<z< 1\\ z, &\text{if }\; 1\leq z< \frac{1}{\theta}\\ \frac{1}{\theta}, &\text{if }\; \frac{1}{\theta} \leq z< \infty \end{cases}
\end{align*}
are increasing in $z\in [0,\infty)$ and $\psi(z)=\frac{\psi_1(z)}{\psi_2(z)}=\frac{\alpha_1+\alpha_2+2}{\alpha_1+\alpha_2+1}\,\frac{1}{\min\big\{1,\frac{1}{z}\big\}}$ is also increasing in $z\in [0,\infty)$. For any fixed $\lambda\geq 1$, 
\begin{align*}
	k_1(z,\lambda)=\frac{\psi\left(\frac{c_{0,1}z}{\nu \lambda}\right)-c_{0,1}}{1-z}
	=\begin{cases} \frac{-\alpha_2}{(\alpha_1+1)(\alpha_1+\alpha_2+1)(1-z)}, &\text{if  } 0<z<\frac{\nu \lambda}{c_{0,1}}\\~\\ \frac{c_{0,1} \left(\frac{\alpha_1+\alpha_2+2}{\alpha_1+\alpha_2+1}\frac{z}{\lambda\,\nu}-1\right)}{1-z},& \text{if  } \frac{\nu \lambda}{c_{0,1}}<z<\infty \end{cases}
\end{align*}
is decreasing in $z\in(0,\infty)$, provided $\nu\geq \frac{\alpha_1+\alpha_2+2}{\alpha_1+\alpha_2+1}$. Since $c_{0,2}=\frac{\alpha_2+2}{\alpha_2+1}>\frac{\alpha_1+\alpha_2+2}{\alpha_1+\alpha_2+1}$, a suitable choice of $\nu$ in this case is $\nu=\frac{\alpha_2+2}{\alpha_2+1}$. To avoid some tedious algebra, we only consider the case $\alpha_1\geq \alpha_2$. For $\alpha_1\geq \alpha_2$ and $\nu=\frac{\alpha_2+2}{\alpha_2+1}$, it is easy to verify that
$$\alpha_1(\lambda)=\alpha_{1,1}=\alpha_{1,\infty}=1-\frac{\alpha_2 (\alpha_2+2)}{2(\alpha_1+2)(\alpha_1+\alpha_2+1)}=\alpha_0,\;\; (say).$$

Using Theorem 4.1.1 (b), we conclude that among the isotonic regression estimators in the class  $\mathcal{D}_{1,\nu}=\{\delta_{1,c_{0,1},\alpha}^{(\nu)}: -\infty <\alpha < \infty\}$ the only admissible estimator is
\begin{align*}
	\delta_{1,c_{0,1},\alpha_0}^{(\nu)}(\underline{X})=\begin{cases}
		\frac{\alpha_1+2}{\alpha_1+1}X_1, &\text{ if  } \frac{\alpha_1+2}{\alpha_1+1}X_1\leq \frac{\alpha_2+2}{\alpha_2+1}X_2\\
		\alpha_0\,\frac{\alpha_1+2}{\alpha_1+1}X_1+(1-\alpha_0)\frac{\alpha_2+2}{\alpha_2+1}X_2, &\text{ if  } \frac{\alpha_1+2}{\alpha_1+1}X_1> \frac{\alpha_2+2}{\alpha_2+1}X_2
	\end{cases};
\end{align*}
here $\alpha_0=1-\frac{\alpha_2 (\alpha_2+2)}{2(\alpha_1+2)(\alpha_1+\alpha_2+1)}$, $\nu=\frac{\alpha_2+2}{\alpha_2+1}$ and $c_{0,1}=\frac{\alpha_1+2}{\alpha_1+1}$. All other estimators in the class $\mathcal{D}_{1,\nu}$, including the BSEE $\delta_{1,c_{0,1}}(\underline{X})=\delta_{1,c_{0,1},1}^{(\nu)}(\underline{X})=\frac{\alpha_1+2}{\alpha_1+1}X_1$, are inadmissible for estimating $\theta_1$ and are dominated by the isotonic regression estimator $\delta_{1,c_{0,1},\alpha_0}^{(\nu)}(\underline{X})$, defined above.

\subsection{\textbf{Isotonic Regression Estimators of Scale Parameter $\theta_2$}}
\label{sec:4.2}
\noindent
\vspace*{2mm}

Let $\beta$ be a fixed real constant, to be suitably chosen as described earlier. Under the notations of Section 4.1, consider the class $\mathcal{D}_{2,\beta}=\{\delta_{2,c_{0,2},\alpha}^{(\beta)}: -\infty <\alpha < \infty\}$ of isotonic regression estimators defined by \eqref{eq:1.8} and \eqref{eq:1.10}; here $\delta_{2,c_{0,2}}(\underline{X})=c_{0,2}X_2$ and $\delta_{1,\beta}(\underline{X})=\beta X_1$. For $z\in S_Z$, let $\psi_3(z)=E[Z_2|Z=z]=z\psi_1(z),$ $\psi_4(z)=E[Z_2^{^2}|Z=z]=z^2\psi_2(z)$ and $\psi^{*}(z)=\frac{\psi_3(z)}{\psi_4(z)}$. For any fixed $\underline{\theta}\in\Theta_0$ (or $\lambda\geq 1$), the risk function
\begin{align*}
	R_2(\underline{\theta},\delta_{2,c_{0,2},\alpha}^{(\beta)})&=E_{\underline{\theta}}\left[\left(\frac{\delta_{2,c_{0,2},\alpha}^{(\beta)}(\underline{X})}{\theta_2}-1\right)^2\right]\\
	&=E_{\underline{\theta}}[(c_{0,2}Z_2-1)^2 I_{(\frac{\beta}{\lambda c_{0,2}},\infty)}(Z)] \\
	&\qquad+ E_{\underline{\theta}}\left[\left(\frac{\alpha\beta}{\lambda}Z_1+(1-\alpha)c_{0,2}Z_2-1\right)^2 I_{(0,\frac{\beta}{\lambda c_{0,2}})}(Z)\right]
\end{align*}
is minimized at $\alpha= \alpha_2(\lambda)$, where
\begin{align}\label{eq:4.3}
	\alpha_2(\lambda)
	&=\frac{1}{c_{0,2}}\frac{E\left[\left(\frac{\beta}{\lambda \,c_{0,2} Z}-1\right)(\psi^{*}(Z)-c_{0,2})\, \psi_4(Z)I_{(0,\frac{\beta}{\lambda c_{0,2}})}(Z)\right]}{E\left[\left(\frac{\beta}{\lambda \,c_{0,2} Z}-1\right)^2\psi_4(Z) I_{(0,\frac{\beta}{\lambda c_{0,2}})}(Z)\right]} \nonumber\\
	=& \frac{1}{c_{0,2}} \frac{\int_{0}^{1} \left(\frac{1}{t}-1\right) (\psi^{*}\left(\frac{\beta t}{\lambda c_{0,2}}\right)-c_{0,2})\, \psi_4\left(\frac{\beta t}{\lambda c_{0,2}}\right)\,f_Z\left(\frac{\beta t}{\lambda c_{0,2}}\right) dt}{\int_{0}^{1} \left(\frac{1}{t}-1\right)^2 \,\psi_4\left(\frac{\beta t}{\lambda c_{0,2}}\right)\,f_Z\left(\frac{\beta t}{\lambda c_{0,2}}\right) dt},\;\;\lambda \geq 1,  \\
	=& \frac{1}{c_{0,2}}E_{\lambda}[k_2(S_{2,\lambda},\lambda)],\;\; \lambda \geq 1,\nonumber
\end{align}
where $k_2(z,\lambda)=\frac{\psi^{*}\left(\frac{\beta z}{\lambda \,c_{0,2}}\right)-c_{0,2}}{\frac{1}{z}-1},\; 0<z<1,\; \lambda \geq 1,$
and $S_{2,\lambda}$ is a r.v. having the pdf
\begin{equation*}
	h_{2,\lambda}(t  )=\begin{cases}
		\frac{\left(\frac{1}{t}-1\right)^2 \,\psi_4\left(\frac{\beta t}{\lambda c_{0,2}}\right)\,f_Z\left(\frac{\beta t}{\lambda c_{0,2}}\right) }{\int_{0}^{1}\left(\frac{1}{s}-1\right)^2 \,\psi_4\left(\frac{\beta s}{\lambda c_{0,2}}\right)\,f_Z\left(\frac{\beta s}{\lambda c_{0,2}}\right) ds},& \text{if}\;\; 0<t<1\\
		0, &\text{otherwise}
	\end{cases},\;\;\lambda \geq 1.
\end{equation*}

It is easy to verify that if, for every fixed $\theta\in(0,1)$, $\frac{\psi_4(\theta t)}{\psi_4(t)}$ and $\frac{f_Z(\theta t)}{f_Z(t)}$ are increasing (decreasing) functions of t on $(0,\frac{\beta}{c_{0,2}})$, then $S_{2,\lambda_1}\leq_{lr}S_{2,\lambda_2}$ ($S_{2,\lambda_2}\leq_{lr}S_{2,\lambda_1}$) and, consequently, $S_{2,\lambda_1}\leq_{st}S_{2,\lambda_2}$ ($S_{2,\lambda_2}\leq_{st}S_{2,\lambda_1}$), whenever $1\leq \lambda_1<\lambda_2<\infty$. \vspace*{2mm}

\noindent
On the lines of Lemma 3.1.1, we have the following result.
\\~\\ \textbf{Lemma 4.2.1 (a)} Suppose that, for every fixed $\theta\in(0,1)$, $\frac{\psi_4(\theta t)}{\psi_4(t)}$ and $\frac{f_Z(\theta t)}{f_Z(t)}$ are increasing (decreasing) functions of $t$ on $(0,\frac{\beta}{c_{0,2}})$. Also, suppose that $\psi^{*}(z)$ is decreasing on $(0,\frac{\beta}{c_{0,2}})$ and, for every fixed $\lambda \geq 1$, $k_2(z,\lambda)=\frac{\psi^{*}\left(\frac{\beta z}{\lambda \,c_{0,2}}\right)-c_{0,2}}{\frac{1}{z}-1}$ is increasing (decreasing) in $z\in(0,1)$. Then $\alpha_2(\lambda)$, defined by (4.3), is an increasing function of $\lambda\in[1,\infty)$,
$\inf_{\lambda \geq 1} \alpha_2(\lambda) = \alpha_2(1)=\alpha_{2,1}, \text{ say, and }
\sup_{\lambda \geq 1} \alpha_2(\lambda) =\lim_{\lambda \to \infty}\alpha_2(\lambda)=\alpha_{2,\infty},\text{ say}.$
\\~\\\textbf{(b)} Suppose that, for every fixed $\theta\in(0,1)$, $\frac{\psi_4(\theta t)}{\psi_4(t)}$ and $\frac{f_Z(\theta t)}{f_Z(t)}$ are increasing (decreasing) functions of $t$ on $(0,\frac{\beta}{c_{0,2}})$. Also, suppose that $\psi^{*}(z)$ is increasing on $(0,\frac{\beta}{c_{0,2}})$ and, for every fixed $\lambda \geq 1$, $k_2(z,\lambda)=\frac{\psi^{*}\left(\frac{\beta z}{\lambda \,c_{0,2}}\right)-c_{0,2}}{\frac{1}{z}-1}$ is decreasing (increasing) in $z\in(0,1)$. Then $\alpha_2(\lambda)$ is a decreasing function of $\lambda\in[1,\infty)$,
$\inf_{\lambda \geq 1} \alpha_2(\lambda)= \lim_{\lambda \to \infty}\alpha_2(\lambda)=\alpha_{2,\infty}, \text{ say, and }
\sup_{\lambda \geq 1} \alpha_2(\lambda) = \alpha_2(1)=\alpha_{2,1} \text{, say}.$
\vspace*{2mm}

\noindent
Now we present the main result of this subsection, whose proof is similar to the proof of Theorem 3.1.1.
\\~\\ \textbf{Theorem 4.2.1 (a)} Suppose that assumptions of Lemma 4.2.1 (a) hold. Then the estimators that are admissible in the class $\mathcal{D}_{2,\beta}$ are $\{\delta_{2,c_{0,2},\alpha}^{(\beta)}:\alpha\in[\alpha_{2,1},\alpha_{2,\infty}]\}$. Moreover, for $-\infty<\alpha_1<\alpha_2\leq \alpha_{2,1}$ or $\alpha_{2,\infty}\leq \alpha_2<\alpha_1$, the estimator $\delta_{2,c_{0,2},\alpha_2}^{(\beta)}$ dominates the estimator $\delta_{2,c_{0,2},\alpha_1}^{(\beta)}$, for any $\underline{\theta}\in\Theta_0$.
\\~\\ \textbf{(b)} Suppose that assumptions of Lemma 4.2.1 (b) hold. Then the estimators that are admissible in the class $\mathcal{D}_{2,\beta}$ are $\{\delta_{2,c_{0,2},\alpha}^{(\beta)}:\alpha\in[\alpha_{2,\infty},\alpha_{2,1}]\}$. Moreover, for $-\infty<\alpha_1<\alpha_2\leq \alpha_{2,\infty}$ or $\alpha_{2,1}\leq \alpha_2<\alpha_1<\infty$, the estimator $\delta_{2,c_{0,2},\alpha_2}^{(\beta)}$ dominates the estimator $\delta_{2,c_{0,2},\alpha_1}^{(\beta)}$, for any $\underline{\theta}\in\Theta_0$.
\\~\\Now we illustrate some applications of Theorem 4.2.1.	
\\~\\ \textbf{Example 4.2.1.} Let $(X_1,X_2)$ be a random vector as defined in Example 4.1.1. Here $c_{0,i}=\frac{E[Z_i]}{E[Z_i^2]}=\frac{1}{\alpha_i+1},\;i=1,2,$ the BSEE of $\theta_2$ is $\delta_{2,0}(\underline{X})=\delta_{2,c_{0,2},0}^{(\beta)}(\underline{X})=\frac{X_2}{\alpha_2+1}$,
$\psi_3(z)=E[Z_2|Z=z]=z\,E[Z_1|Z=z]=(\alpha_1+\alpha_2)\frac{z}{1+z},\;z>0,\;
\psi_4(z)=E[Z_2^{^2}|Z=z]=z^2\,E[Z_1^{^2}|Z=z]=(\alpha_1+\alpha_2+1)(\alpha_1+\alpha_2)\frac{z^2}{(1+z)^2},\; z>0,$ and $\psi^{*}(z)=\frac{\psi_3(z)}{\psi_4(z)}=\frac{1+z}{(\alpha_1+\alpha_2+1)z}$.
Clearly, for $0<\theta<1$, $\frac{f_Z(\theta z)}{f_Z(z)}$ and $\frac{\psi_2(\theta z)}{\psi_2(z)}=\frac{\theta^2(1+z)^2}{(1+\theta z)^2}$ are increasing in $z\in \left[0,\frac{\beta}{c_{0,2}}\right]$ and $\psi^{*}(z)=\frac{\psi_3(z)}{\psi_4(z)}=\frac{1+z}{(\alpha_1+\alpha_2+1)z}$ is decreasing in $z\in \left[0,\frac{\beta}{c_{0,2}}\right]$. For any fixed $\lambda\geq 1$,
$ k_2(z,\lambda)=\frac{\psi^{*}\left(\frac{\beta z}{\lambda\, c_{0,2}}\right)-c_{0,2}}{\frac{1}{z}-1}=\frac{\lambda-\alpha_1 \beta z}{ (\alpha_1+\alpha_2+1)(\alpha_2+1)\beta (1-z)},\; z>0$
is increasing in $z\in[0,1]$, provided $\beta\leq \frac{1}{\alpha_1}$. Since $c_{0,1}=\frac{1}{\alpha_1+1}$, a suitable choice of $\beta$ in this case is $\beta=\frac{1}{\alpha_1+1}$. For $\beta=\frac{1}{\alpha_1+1}$, it is easy to verify that
\begin{align*}
	\alpha_{2,1}&=\alpha_2(1)\\
	&=\frac{1}{c_{0,2}}\frac{\int_{0}^{1} \left(\frac{1}{t}-1\right) \left(\psi^{*}\left(\frac{\beta t}{ c_{0,2}}\right)-c_{0,2}\right) \psi_4\left(\frac{\beta t}{c_{0,2}}\right) f_Z\left(\frac{\beta t}{c_{0,2}}\right) dt}{\int_{0}^{1}  \left(\frac{1}{t}-1\right)^2\, \psi_4\left(\frac{\beta t}{c_{0,2}}\right) f_Z\left(\frac{\beta t}{c_{0,2}}\right) dt}\\
	&=\frac{1}{\alpha_1+\alpha_2+1}\left(\alpha_1+1+
	\frac{\int_{0}^{1} \left(\frac{1}{t}-1\right) \frac{1}{\left(1+\frac{\alpha_1+1}{(\alpha_2+1) t}\right)^{\alpha_1+\alpha_2+2}}\frac{1}{t^{\alpha_1+1}} dt}{\int_{0}^{1} \left(\frac{1}{t}-1\right)^2 \frac{1}{\left(1+\frac{\alpha_1+1}{(\alpha_2+1) t}\right)^{\alpha_1+\alpha_2+2}}\frac{1}{t^{\alpha_1+1}} dt}\right)
\end{align*}
\begin{align*}	\text{and}\qquad
	\alpha_{2,\infty}=\lim_{\lambda \to \infty} \alpha_2(\lambda)=\infty.\qquad\qquad\qquad\qquad\qquad\qquad\qquad\qquad
\end{align*}
Using Theorem 4.2.1. (b), it follows that the estimators $\{\delta_{2,c_{0,2},\alpha}^{(\beta)}: \alpha\in[\alpha_{2,1},\infty)\}$ are admissible within the class  $\mathcal{D}_{2,\beta}=\{\delta_{2,c_{0,2},\alpha}^{(\beta)}: -\infty <\alpha < \infty\}$ of isotonic regression estimators of $\theta_2$; here $\beta=\frac{1}{\alpha_1+1}$ and $c_{0,2}=\frac{1}{\alpha_2+1}$.  Moreover, for $-\infty<\alpha_1<\alpha_2\leq \alpha_{2,1}$, the estimator $\delta_{2,c_{0,2},\alpha_2}^{(\beta)}$ dominates the estimator $\delta_{2,c_{0,2},\alpha_1}^{(\beta)}$. Since $\alpha_{2,1}>\frac{\alpha_1+1}{\alpha_1+\alpha_2+1}>\frac{\alpha_1}{\alpha_1+\alpha_2}$, it follows that the BSEE $\delta_{2,c_{0,2}}(\underline{X})=\delta_{2,c_{0,2},0}^{(\beta)}(\underline{X})=\frac{X_2}{\alpha_2+1}$ is inadmissible for estimating $\theta_2$ and is dominated by
\begin{align*}
	\delta_{2,c_{0,2},\beta_0}^{(\beta)}(\underline{X})&=\begin{cases}
		\frac{X_2}{\alpha_2+1}, &\text{ if  } \frac{X_1}{\alpha_1+1}\leq \frac{X_2}{\alpha_2+1}\\
		\alpha_{2,1}\frac{X_1}{\alpha_1+1}+(1-\alpha_{2,1})\frac{X_2}{\alpha_2+1}, &\text{ if  } \frac{X_1}{\alpha_1+1}> \frac{X_2}{\alpha_2+1}
	\end{cases}\\
	&=\max\bigg\{\frac{X_2}{\alpha_2+1},\alpha_{2,1}\frac{X_1}{\alpha_1+1}+(1-\alpha_{2,1})\frac{X_2}{\alpha_2+1}\bigg\}.
\end{align*}
\textbf{Example 4.2.2.} Let $X_1$ and $X_2$ be independent random variables, as defined in Example 4.1.2. Here $c_{0,i}=\frac{E[Z_i]}{E[Z_i^2]}=\frac{\alpha_i+2}{\alpha_i+1},\;i=1,2,$ and, for any fixed $\beta>0$, the BSEE of $\theta_2$ is $\delta_{2,0}(\underline{X})=\delta_{2,c_{0,2},0}^{(\beta)}(\underline{X})=\frac{\alpha_2+2}{\alpha_2+1}X_2$. $\text{Moreover, }
\psi_3(z)=E[Z_2|Z=z]=zE[Z_1|Z=z]= \frac{\alpha_1+\alpha_2}{\alpha_1+\alpha_2+1}\min\{z,1\},\, z>0,\;
\psi_4(z)=E[Z_2^2|Z=z]=z^2 E[Z_1^2|Z=z]=\frac{\alpha_1+\alpha_2}{\alpha_1+\alpha_2+2}\left(\min\{z,1\}\right)^2,\, z>0,
\text{ and, for }0<\theta<1$, 
$$\frac{\psi_4(\theta z)}{\psi_4(z)}=\frac{\left(\min\{\theta z,1\}\right)^2}{\left(\min\{z,1\}\right)^2}=\begin{cases} \theta^2, &\text{if } \; 0<z< 1\\ \theta^2 z^2, &\text{if }\; 1\leq z< \frac{1}{\theta}\\ 1, &\text{if }\; \frac{1}{\theta} \leq z< \infty \end{cases}$$
is increasing in $z\in \left[0,\frac{\beta}{c_{0,2}}\right]$, where $\beta>0$ is fixed. Also, $\psi^{*}(z)=\frac{\psi_3(z)}{\psi_4(z)}=\frac{\alpha_1+\alpha_2+2}{\alpha_1+\alpha_2+1}\,\frac{1}{\min\{z,1\}}$ is decreasing in $z\in \left[0,\frac{\beta}{c_{0,2}}\right]$. For any fixed $\lambda\geq 1$,

\begin{align*}
	k_2(z,\lambda)&=\frac{\psi\left(\frac{\beta z}{c_{0,2}\, \lambda}\right)-c_{0,2}}{\frac{1}{z}-1}
	=\begin{cases}\frac{\alpha_2+2}{\alpha_2+1}\frac{ \left(\frac{\alpha_1+\alpha_2+2}{\alpha_1+\alpha_2+1}\frac{\lambda}{\beta}-z\right)}{1-z}, &\text{if  } 0<z<\frac{\lambda \,c_{0,2}}{\beta}\\~\\ \frac{-\alpha_1}{(\alpha_2+1)(\alpha_1+\alpha_2+1)\left(\frac{1}{z}-1\right)},& \text{if  } \frac{\lambda\,c_{0,2}}{\beta}<z<\infty \end{cases}
\end{align*}

is increasing in $z\in[0,1]$, provided $\beta\leq \frac{\alpha_1+\alpha_2+2}{\alpha_1+\alpha_2+1}$. Since $c_{0,1}=\frac{\alpha_1+2}{\alpha_1+1}>\frac{\alpha_1+\alpha_2+2}{\alpha_1+\alpha_2+1}$, an appropriate choice of $\beta$ in this case is $\beta=\frac{\alpha_1+\alpha_2+2}{\alpha_1+\alpha_2+1}$. For $\beta=\frac{\alpha_1+\alpha_2+2}{\alpha_1+\alpha_2+1}$, it is easy to verify that
$$	\alpha_2(\lambda)=\frac{\int_{0}^{1}(1-t) (\lambda-t) t^{\alpha_2-1}dt}{\int_{0}^{1}(1-t)^2\, t^{\alpha_2-1}dt},\;\lambda\geq 0,$$
$$	\alpha_{2,1}=\alpha_2(1)=1\qquad \text{and} \qquad \alpha_{2,\infty}=\lim_{\lambda \to \infty}\alpha_2(\lambda)=\infty.
$$
Using Theorem 4.2.1. (a), it follows that the estimators $\{\delta_{2,c_{0,2},\alpha}^{(\beta)}: \alpha\in [1,\infty)\}$ are admissible within the class  $\mathcal{D}_{2,\beta}=\{\delta_{2,c_{0,2},\alpha}^{(\beta)}: -\infty <\alpha < \infty\}$ of isotonic regression estimators of $\theta_2$; here $\beta=\frac{\alpha_1+\alpha_2+2}{\alpha_1+\alpha_2+1}$ and $c_{0,2}=\frac{\alpha_2+2}{\alpha_2+1}$. Moreover, for $-\infty<\alpha_1<\alpha_2\leq 1$, the estimator $\delta_{2,c_{0,2},\alpha_2}^{(\beta)}$ dominates the estimator $\delta_{2,c_{0,2},\alpha_1}^{(\beta)}$. In particular the BSEE $\delta_{2,c_{0,2}}(\underline{X})=\delta_{2,c_{0,2},0}^{(\beta)}(\underline{X})=\frac{\alpha_2+2}{\alpha_2+1}X_2$ is inadmissible for estimating $\theta_2$ and is dominated by the isotonic regression estimator
\begin{align*}
	\delta_{2,c_{0,2},1}^{(\beta)}(\underline{X})&=\begin{cases}
		\frac{\alpha_2+2}{\alpha_2+1}X_2, &\text{ if  } \frac{\alpha_1+\alpha_2+2}{\alpha_1+\alpha_2+1}X_1\leq \frac{\alpha_2+2}{\alpha_2+1}X_2\\
		\frac{\alpha_1+\alpha_2+2}{\alpha_1+\alpha_2+1}X_1, &\text{ if  } \frac{\alpha_1+\alpha_2+2}{\alpha_1+\alpha_2+1}X_1> \frac{\alpha_2+2}{\alpha_2+1}X_2
	\end{cases}
	=\max\bigg\{\frac{\alpha_2+2}{\alpha_2+1}X_2,\frac{\alpha_1+\alpha_2+2}{\alpha_1+\alpha_2+1}X_1\bigg\}.
\end{align*}
This is another situation where the isotonic regression estimator based on the BSEE $(\frac{\alpha_1+2}{\alpha_1+1}X_1,\frac{\alpha_2+2}{\alpha_2+1}X_2)$ may not be able to provide improvement over the BSEE $\delta_{2,c_{0,2}}(\underline{X})=\frac{\alpha_2+2}{\alpha_2+1}X_2$. Rather an isotonic regression estimator based on $(\frac{\alpha_1+\alpha_2+2}{\alpha_1+\alpha_2+1}X_1,\frac{\alpha_2+2}{\alpha_2+1}X_2)$ provides improvement over the BSEE $\delta_{2,c_{0,2}}(\underline{X})$.

\subsection{\textbf{Simulation Study For Estimation of Scale Parameter $\theta_1$}} \label{sec:4.3}
\noindent
\vspace*{2mm}

\begin{figure}[h!]
	\begin{subfigure}{0.48\textwidth}
		\centering
		\includegraphics[width=72mm,scale=1.2]{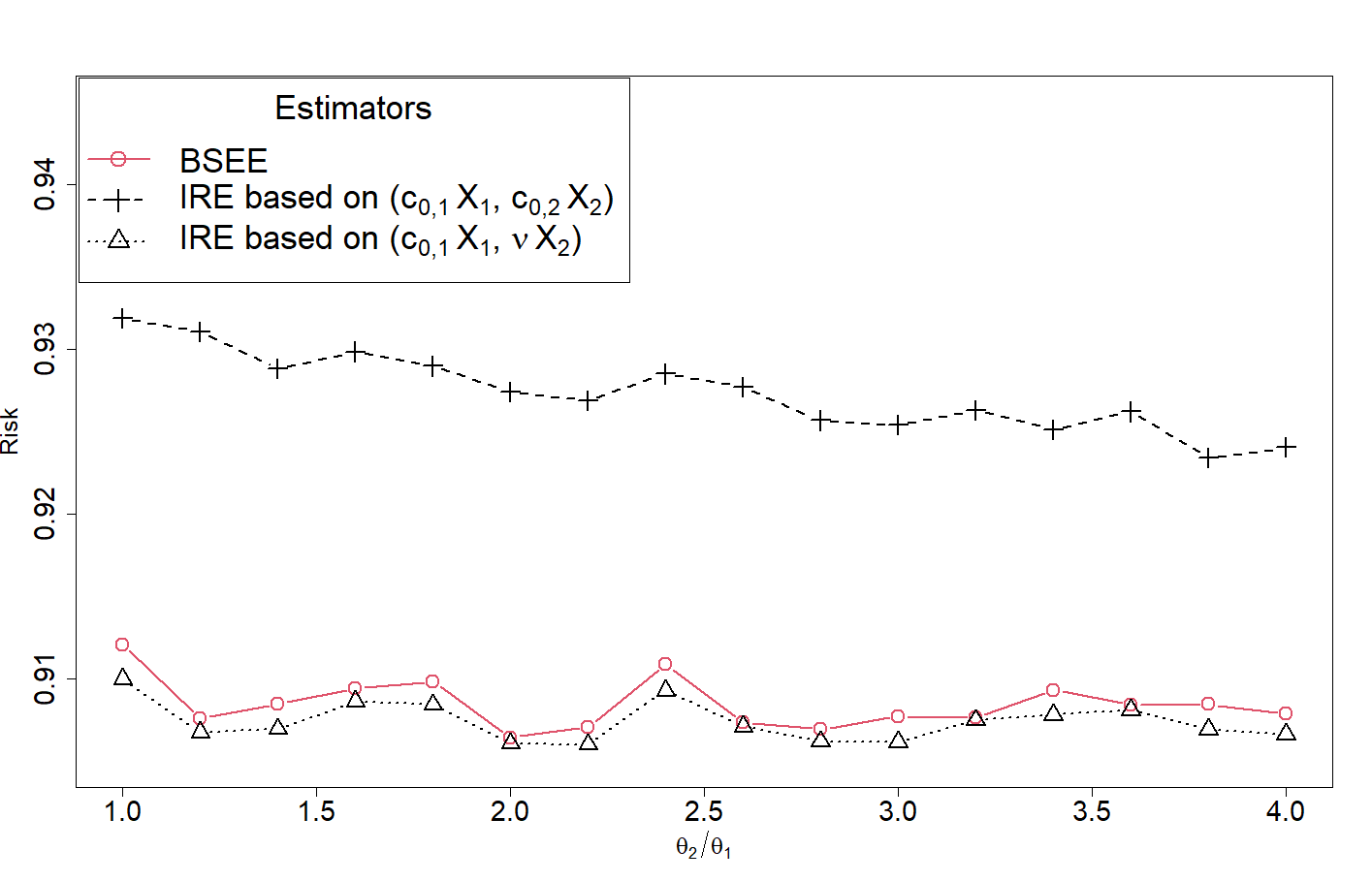} 
		\caption{ $\alpha_1=0.1$ and $\alpha_2=0.2$.} 
		\label{fig2:a}  
	\end{subfigure}
	\begin{subfigure}{0.48\textwidth}
		\centering
		
		\includegraphics[width=72mm,scale=1.2]{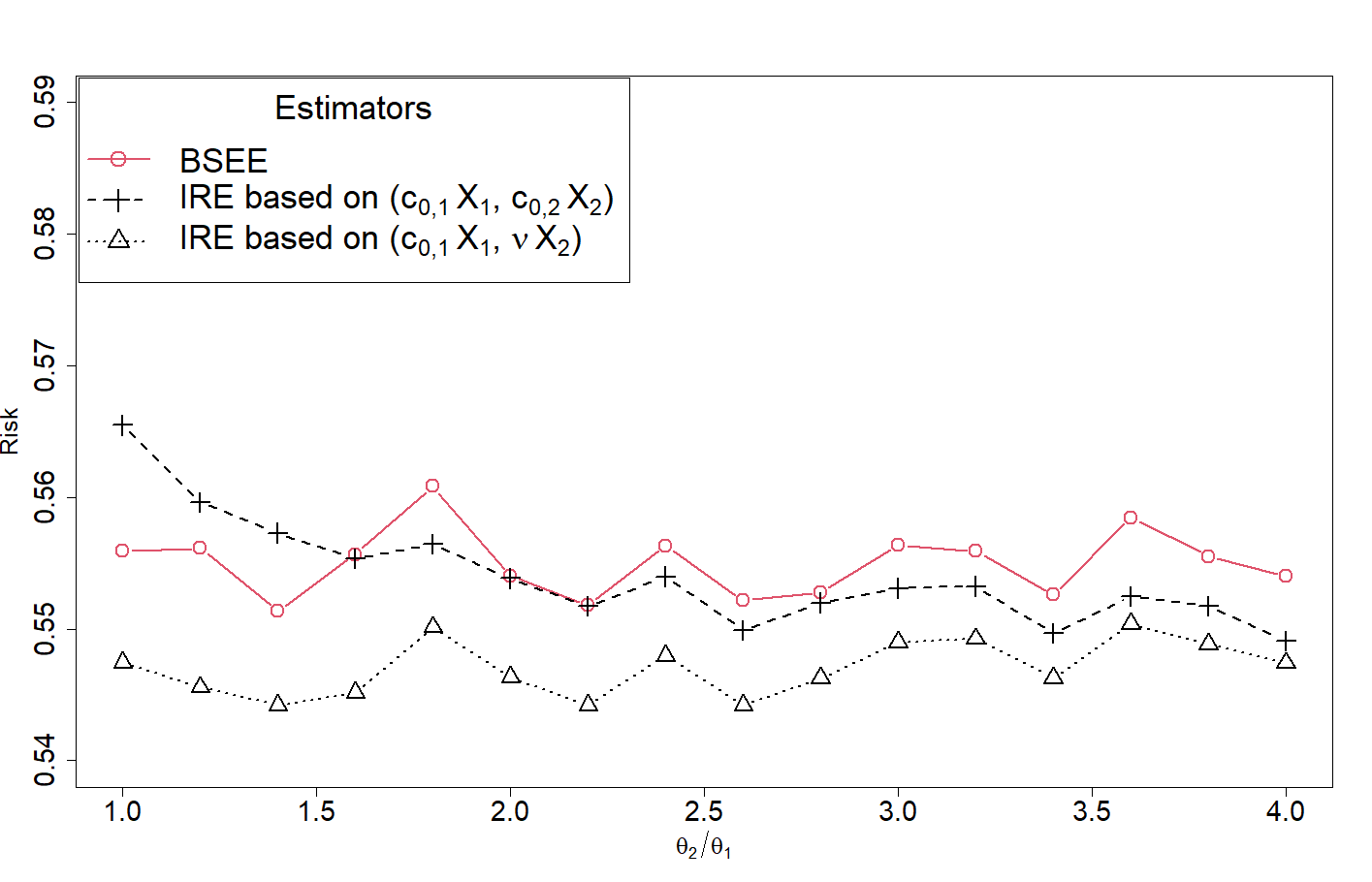} 
		
		\caption{ $\alpha_1=0.8$ and $\alpha_2=0.5$.} 
		\label{fig2:b} 
	\end{subfigure}
	\\	\begin{subfigure}{0.48\textwidth}
		\centering
		\includegraphics[width=72mm,scale=1.2]{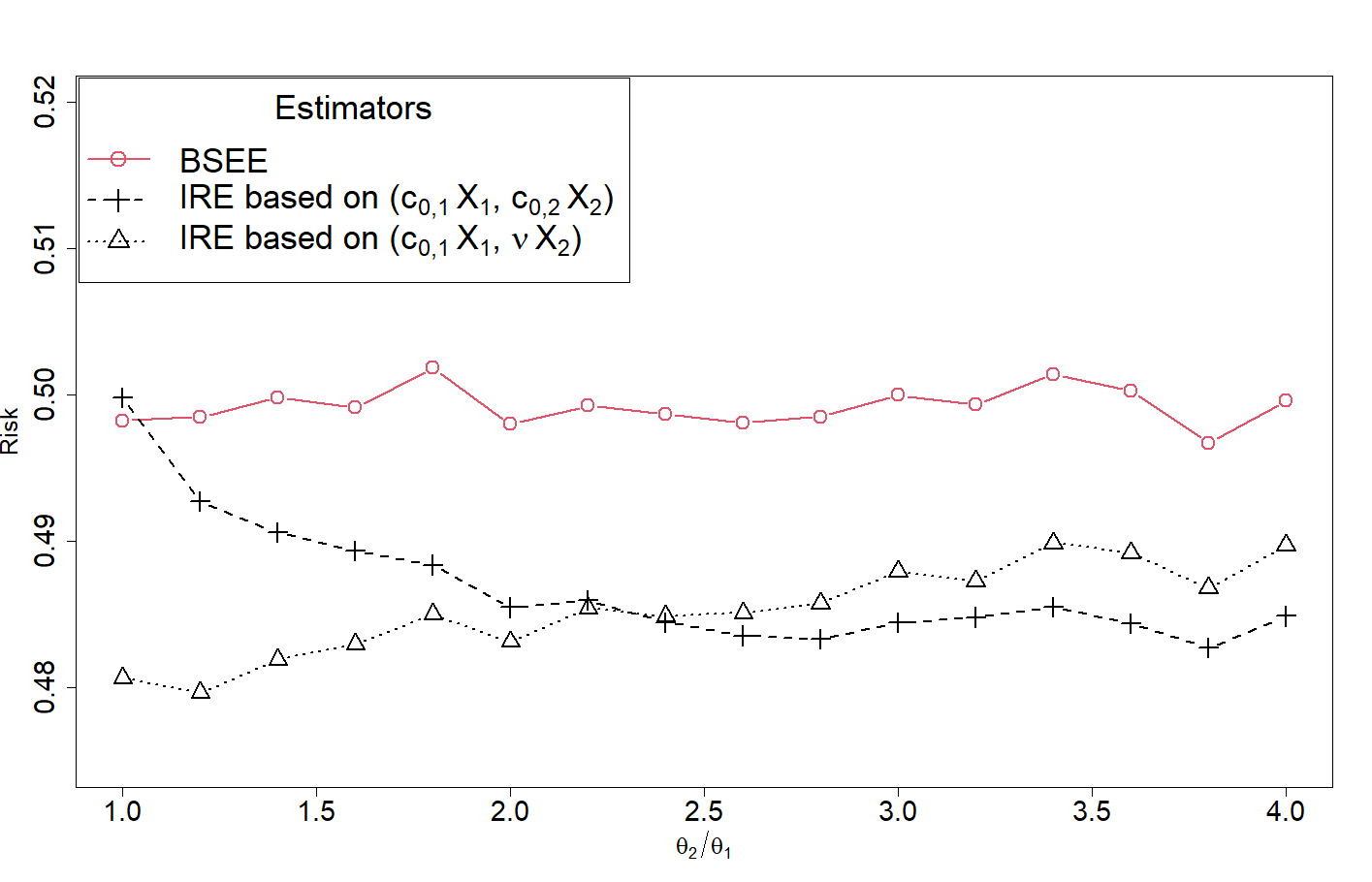} 
		\caption{ $\alpha_1=1$ and $\alpha_2=1$.} 
		\label{fig2:c} 
	\end{subfigure}
	\begin{subfigure}{0.48\textwidth}
		\centering
		\includegraphics[width=72mm,scale=1.2]{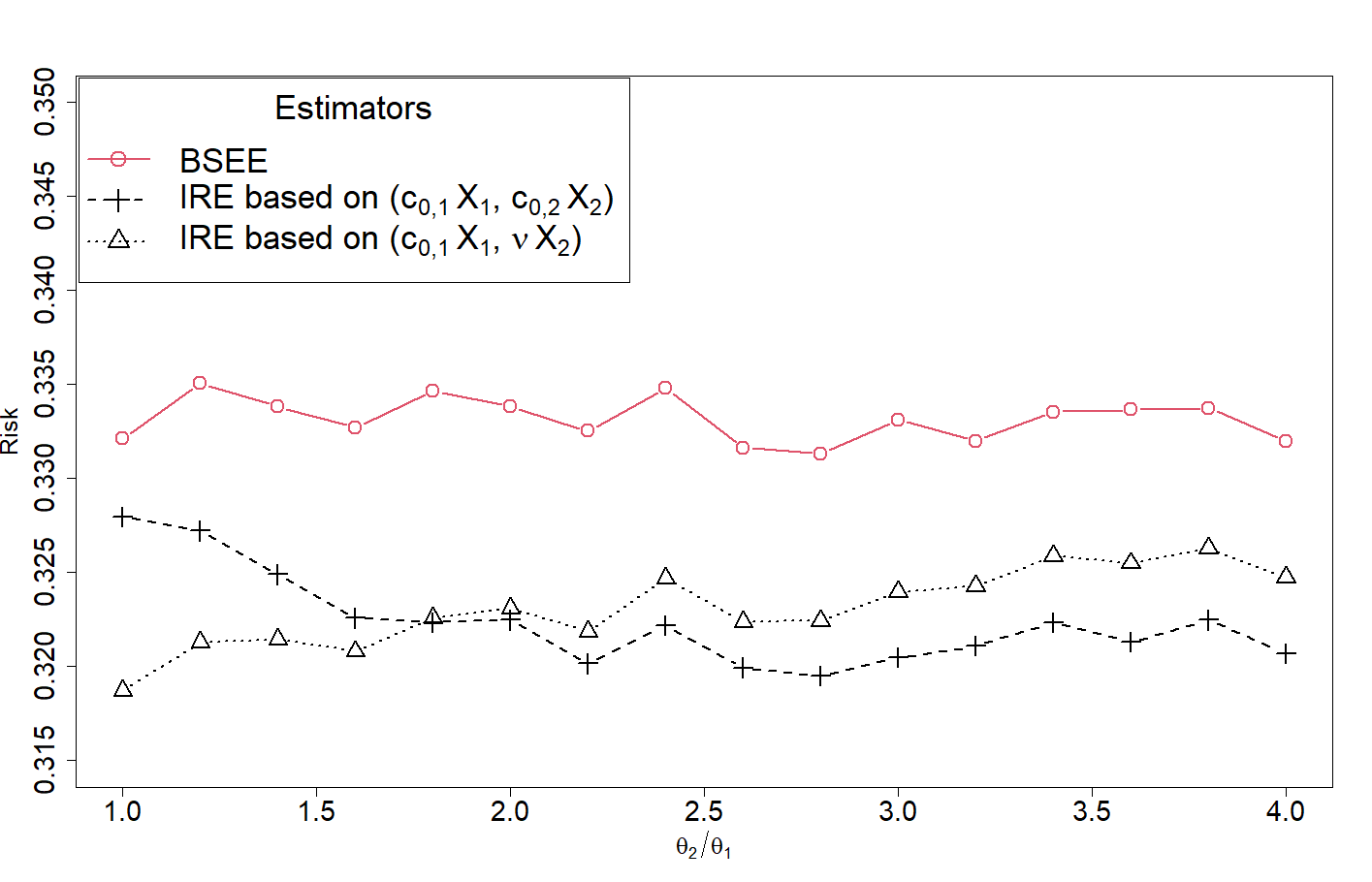} 
		\caption{ $\alpha_1=2$ and $\alpha_2=1$.} 
		\label{fig2:d} 
	\end{subfigure}

	\caption{Risk plots of the BSEE $c_{0,1}X_1$, the IRE based on $(c_{0,1}X_1,c_{0,2}X_2)$ with $\alpha=\frac{\alpha_1}{\alpha_1+\alpha_2}$, and the IRE based on $(c_{0,1}X_1,\nu X_2)$ with $\alpha=\beta_0$, where $\nu=\frac{1}{\alpha_2}$ and $c_{0,i}=\frac{1}{\alpha_i+1},\;i=1,2$, against the values of $\lambda=\frac{\theta_2}{\theta_1}$.}
	\label{fig2}
\end{figure}

In Example 4.1.1, we have considered two independent gamma distributions with unknown order restricted scale parameters (i.e., $\theta_1\leq \theta_2$) and known shape parameters ($\alpha_1>0$ and $\alpha_2>0$). For estimation of smaller scale parameter $\theta_1$, under the scaled squared error loss function, we showed that the IRE $\min\{\frac{X_1}{\alpha_1},\frac{X_1+X_2}{\alpha_1+\alpha_2+1}\}$, based on ($\frac{X_1}{\alpha_1+1},\frac{X_2}{\alpha_2}$) with $\alpha=\beta_0=\frac{\alpha_1+1}{\alpha_1+\alpha_2+1}$, dominates the BLEE $\frac{X_1}{\alpha_1+1}$. To further evaluate the performances of various estimators, under the scaled squared error loss function, in this section, we compare the risk performances of the BSEE $\frac{X_1}{\alpha_1+1}$, the IRE based on the BSEEs ($\frac{X_1}{\alpha_1+1},\frac{X_2}{\alpha_2+1}$) with $\alpha=\frac{\alpha_1}{\alpha_1+\alpha_2}$, and the IRE based on ($\frac{X_1}{\alpha_1+1},\frac{X_2}{\alpha_2}$) with $\alpha=\frac{\alpha_1+1}{\alpha_1+\alpha_2+1}$, numerically, through Monte Carlo simulations. For simulations, we have generated 10000 sample of size 1 each from relevant gamma distributions, for different values of known shape parameters ($\alpha_1$ and $\alpha_2$), and computed the simulated risks of the BSEE, the IRE based on ($\frac{X_1}{\alpha_1+1},\frac{X_2}{\alpha_2+1}$) with $\alpha=\frac{\alpha_1}{\alpha_1+\alpha_2}$, and the IRE based on ($\frac{X_1}{\alpha_1+1},\frac{X_2}{\alpha_2}$) with $\alpha=\frac{\alpha_1+1}{\alpha_1+\alpha_2+1}$.
The simulated values of risks of various estimators are plotted in Figure \ref{fig2}. The following observations are evident from Figure \ref{fig2}:
\\~\\(i) The risk function values of the IRE based on ($\frac{X_1}{\alpha_1+1},\frac{X_2}{\alpha_2}$) with $\alpha=\frac{\alpha_1+1}{\alpha_1+\alpha_2+1}$ is nowhere larger than the risk function values of the BSEE $\frac{X_1}{\alpha+1}$, which is in conformity with theoretical findings of Example 4.1.3.
\\(ii) As we mentioned in the Example 4.1.1, here isotonic regression estimators (IREs) based on the BSEEs ($\frac{X_1}{\alpha_1+1},\frac{X_2}{\alpha_2+1}$) does not always provide improvements over the BSEE $\frac{X_1}{\alpha_1+1}$.
\noindent
\\(iii) There is no clear cut winner between the the IRE based on ($\frac{X_1}{\alpha_1+1},\frac{X_2}{\alpha_2+1}$) with $\alpha=\frac{\alpha_1}{\alpha_1+\alpha_2}$ and the IRE based on ($\frac{X_1}{\alpha_1+1},\frac{X_2}{\alpha_2}$) with $\alpha=\frac{\alpha_1+1}{\alpha_1+\alpha_2+1}$. But the IRE based on ($\frac{X_1}{\alpha_1+1},\frac{X_2}{\alpha_2}$) with $\alpha=\frac{\alpha_1+1}{\alpha_1+\alpha_2+1}$ performs reasonably well in comparison to other estimators. 
\FloatBarrier
\FloatBarrier

\section*{\textbf{Funding}}

This work was supported by the [Council of Scientific and Industrial Research (CSIR)] under Grant [number 09/092(0986)/2018].

\bibliographystyle{apalike}
\bibliography{Paper3}	

\begin{thebibliography}{}

\bibitem[Artin, 1931]{artin}
Artin, E. (1931).
\newblock Einf{\"u}hrung in die theorie der gammafunktion, hamburger math.

\bibitem[Barlow et~al., 1972]{MR0326887}
Barlow, R.~E., Bartholomew, D.~J., Bremner, J.~M., and Brunk, H.~D. (1972).
\newblock {\em Statistical inference under order restrictions. {T}he theory and
  application of isotonic regression}.
\newblock John Wiley \& Sons.

\bibitem[Chang et~al., 2017]{MR3698503}
Chang, Y.-T., Fukuda, K., and Shinozaki, N. (2017).
\newblock Estimation of two ordered normal means when a covariance matrix is
  known.
\newblock {\em Statistics}, 51(5):1095--1104.

\bibitem[Hwang and Peddada, 1994]{MR1272076}
Hwang, J. T.~G. and Peddada, S.~D. (1994).
\newblock Confidence interval estimation subject to order restrictions.
\newblock {\em Ann. Statist.}, 22(1):67--93.

\bibitem[Kelly, 1989]{MR994278}
Kelly, R.~E. (1989).
\newblock Stochastic reduction of loss in estimating normal means by isotonic
  regression.
\newblock {\em Ann. Statist.}, 17(2):937--940.

\bibitem[Kubokawa and Saleh, 1994]{MR1370413}
Kubokawa, T. and Saleh, A. K. M.~E. (1994).
\newblock Estimation of location and scale parameters under order restrictions.
\newblock {\em J. Statist. Res.}, 28(1-2):41--51.

\bibitem[Kumar and Sharma, 1988]{MR981031}
Kumar, S. and Sharma, D. (1988).
\newblock Simultaneous estimation of ordered parameters.
\newblock {\em Comm. Statist. Theory Methods}, 17(12):4315--4336.

\bibitem[Kushary and Cohen, 1989]{MR1029476}
Kushary, D. and Cohen, A. (1989).
\newblock Estimating ordered location and scale parameters.
\newblock {\em Statist. Decisions}, 7(3):201--213.

\bibitem[Lee, 1981]{MR615447}
Lee, C. I.~C. (1981).
\newblock The quadratic loss of isotonic regression under normality.
\newblock {\em Ann. Statist.}, 9(3):686--688.

\bibitem[Marshall and Olkin, 1979]{MR552278}
Marshall, A.~W. and Olkin, I. (1979).
\newblock {\em Inequalities: theory of majorization and its applications},
  volume 143 of {\em Mathematics in Science and Engineering}.
\newblock Academic Press, Inc. [Harcourt Brace Jovanovich, Publishers], New
  York-London.

\bibitem[Marshall and Olkin, 2007]{MR2363282}
Marshall, A.~W. and Olkin, I. (2007).
\newblock Characterizations of distributions through coincidences of
  semiparametric families.
\newblock {\em J. Statist. Plann. Inference}, 137(11):3618--3625.

\bibitem[Misra and Singh, 1994]{MR1366828}
Misra, N. and Singh, H. (1994).
\newblock Estimation of ordered location parameters: the exponential
  distribution.
\newblock {\em Statistics}, 25(3):239--249.

\bibitem[Misra and van~der Meulen, 2003]{MR1985890}
Misra, N. and van~der Meulen, E.~C. (2003).
\newblock On stochastic properties of {$m$}-spacings.
\newblock {\em J. Statist. Plann. Inference}, 115(2):683--697.

\bibitem[Pal and Kushary, 1992]{MR1165709}
Pal, N. and Kushary, D. (1992).
\newblock On order restricted location parameters of two exponential
  distributions.
\newblock {\em Statist. Decisions}, 10(1-2):133--152.

\bibitem[Patra and Kumar, 2017]{patra}
Patra, L.~K. and Kumar, S. (2017).
\newblock Estimating ordered means of a bivariate normal distribution.
\newblock {\em American Journal of Mathematical and Management Sciences},
  36(2):118--136.

\bibitem[Pe\v{c}ari\'{c} et~al., 1992]{MR1162312}
Pe\v{c}ari\'{c}, J.~E., Proschan, F., and Tong, Y.~L. (1992).
\newblock {\em Convex functions, partial orderings, and statistical
  applications}, volume 187 of {\em Mathematics in Science and Engineering}.
\newblock Academic Press, Inc., Boston, MA.

\bibitem[Pr\'{e}kopa, 1971]{MR315079}
Pr\'{e}kopa, A. (1971).
\newblock Logarithmic concave measures with application to stochastic
  programming.
\newblock {\em Acta Sci. Math. (Szeged)}, 32:301--316.

\bibitem[Robertson et~al., 1988]{MR961262}
Robertson, T., Wright, F.~T., and Dykstra, R.~L. (1988).
\newblock {\em Order restricted statistical inference}.
\newblock John Wiley \& Sons.

\bibitem[Shaked and Shanthikumar, 2007]{MR2265633}
Shaked, M. and Shanthikumar, J.~G. (2007).
\newblock {\em Stochastic orders}.
\newblock Springer, New York.

\bibitem[van Eeden, 2006]{MR2265239}
van Eeden, C. (2006).
\newblock {\em Restricted parameter space estimation problems. Admissibility
  and minimaxity properties}, volume 188 of {\em Lecture Notes in Statistics}.
\newblock Springer, New York.

\bibitem[Vijayasree et~al., 1995]{MR1345425}
Vijayasree, G., Misra, N., and Singh, H. (1995).
\newblock Componentwise estimation of ordered parameters of {$k\ (\geq 2)$}
  exponential populations.
\newblock {\em Ann. Inst. Statist. Math.}, 47(2):287--307.

\end{thebibliography}
\end{document}